\documentclass[11pt,amssymb,amsmath]{article}
\usepackage{epsf}
\usepackage{amssymb}
\usepackage{graphics}
\usepackage{amsthm}
\usepackage{amsfonts}
\usepackage{latexsym}
\usepackage{mathrsfs} 
\usepackage{graphicx} 
\usepackage{color} 
\usepackage{xcolor} 
\usepackage{amsxtra} 
\usepackage{mathtools}

\begin{document}

\newcommand{\om}{\omega}

\newcommand{\beqn}{\begin{eqnarray}}
\newcommand{\eeqn}{\end{eqnarray}}
\newcommand{\fr}{\frac}
\newcommand{\nin}{\noindent}
\newcommand{\beq}{\begin{equation}}
\newcommand{\eeq}{\end{equation}}
\newtheorem{defn}{Definition}[section]
\newtheorem{thm}[defn]{Theorem}
\newtheorem{prop}[defn]{Proposition}
\newtheorem{cor}[defn]{Corollary}
\newtheorem{remark}[defn]{Remark}
\newtheorem{remarks}[defn]{Remarks}
\newtheorem{conjecture}[defn]{Conjecture}
\newtheorem{example}[defn]{Example}
\newtheorem{problem}{Problem}
\newtheorem{defn-prop}[defn]{Definition--Proposition}

\newenvironment{rem}{\smallskip\noindent%
\refstepcounter{subsection}%
{\bf \thesubsection}~~{\sc Remark.}\hspace{-1mm}}
{\smallskip}

\newenvironment{ex}{\smallskip\noindent%
\refstepcounter{subsection}%
{\bf \thesubsection}~~{\sc Example.}}
{\smallskip}

\newenvironment{que}{\smallskip\noindent%
\refstepcounter{subsection}%
{\bf \thesubsection}~~{\sc Question.}\hspace{-1mm}}
{\smallskip}

\renewcommand{\theequation}{\arabic{section}.\arabic{equation}}

\newcommand{\Vir}{\text{\upshape Vir}}
\newcommand{\Rot}{\text{\upshape Rot}}
\newcommand{\Diff}{\mathscr{D}}
\newcommand{\Vect}{\text{\upshape Vect}}
\newcommand{\Met}{\text{\upshape Met}}
\newcommand{\Metmu}{\Met_{\mu}}
\newcommand{\Vol}{\text{\upshape Dens}}
\newcommand{\SDiff}{\mathscr{D}_\mu}
\newcommand{\SVect}{\text{\upshape SVect}_\mu}
\newcommand{\Ham}{\text{\upshape Ham}}
\newcommand{\HamM}{\Ham_\omega(M)}
\newcommand{\ham}{\text{\upshape ham}}
\newcommand{\hamM}{\ham_\omega(M)}
\newcommand{\Symp}{\text{\upshape Symp}}
\newcommand{\symp}{\text{\upshape symp}}
\newcommand{\SDiffM}{\mathscr{D}_\mu(M)}
\newcommand{\SympM}{\Symp_\omega(M)}
\newcommand{\sympM}{\symp_\omega(M)}
\newcommand{\Diffmu}{\mathscr{D}_\mu}
\newcommand{\Isog}{\text{Iso}_g}
\newcommand{\DiffM}{\Diff(M)}
\newcommand{\MetM}{\Met(M)}
\newcommand{\MetmuM}{\Metmu(M)}
\newcommand{\VolM}{\Vol(M)}
\newcommand{\DiffmuM}{\mathscr{D}_\mu(M)}
\newcommand{\R}{\mathbb R}
\newcommand{\C}{\mathbb C}


\newcommand{\marginnote}[1]
{
}

\newcounter{gm}
\newcommand{\gm}[1]
{\stepcounter{gm}$^{\bf GM\thegm}$%
\footnotetext{\hspace{-3.7mm}$^{\blacksquare\!\blacksquare}$
{\bf GM\thegm:~}#1}}

\newcounter{bk}
\newcommand{\bk}[1]
{\stepcounter{bk}$^{\bf BK\thebk}$%
\footnotetext{\hspace{-3.7mm}$^{\blacksquare\!\blacksquare}$
{\bf BK\thebk:~}#1}}

\newcounter{as}
\newcommand{\as}[1]
{\stepcounter{as}$^{\bf AS\theas}$%
\footnotetext{\hspace{-3.7mm}$^{\blacksquare\!\blacksquare}$
{\bf AS\theas:~}#1}}


\title{Geometric Hydrodynamics in Open Problems}

\author{Boris Khesin,  Gerard Misio\l ek,  and  Alexander Shnirelman}

\date{~} 
\maketitle 

\numberwithin{equation}{section} 

\begin{abstract} 
Geometric Hydrodynamics has flourished ever since the celebrated 1966 paper of V. Arnold. 
In this paper we present a collection of open
problems along with several new constructions in fluid dynamics and a
concise survey of recent developments and achievements in this area.
The topics discussed include variational settings for different types of fluids, models for invariant metrics,
the Cauchy and boundary value problems, partial analyticity of solutions to the Euler equations,
their steady and singular vorticity solutions, 
differential and Hamiltonian geometry of diffeomorphism groups, long-time behaviour of fluids, as well as
mechanical models of direct and inverse cascades.
%
\end{abstract} 

\tableofcontents

\section{Introduction and notations} 
\label{s:intro} 
We start with the basic setup of geometric hydrodynamics. 
It will provide the background and motivation for various developments discussed in the sequel. 

Consider an ideal (i.e., incompressible and inviscid) fluid in a fixed domain $M$ in $\mathbb{R}^n$ 
($n=2, 3$). In the Eulerian representation a fluid motion is described by 
the evolution of its velocity field which satisfies the  {\it incompressible Euler equations} 
\begin{equation}  \label{ideal}
\partial_t v+v\cdot \nabla v=-\nabla p, 
\quad 
\textrm{div}\, v = 0 \,,
\end{equation}
where the velocity field $v=v(t,x)$ is assumed to be tangent to the domain's boundary, if $\partial M \neq \emptyset$. 
The pressure function $p=p(t,x)$ on the right-hand side is defined uniquely 
by these  conditions, modulo an additive constant, 
and 
the divergence of the field $v$ is computed with respect to the volume form $\mu$ in $\mathbb{R}^n$.

While some ideas can be traced back to  Helmholtz and Kelvin, the modern geometric approach 
to hydrodynamics began with the seminal 1966 paper of Arnold \cite{arn66}. 
It is based on the Lagrangian representation of fluid flows in terms of particle trajectories 
which can be viewed as curves in the infinite-dimensional configuration space given by the group $\mathscr{D}_\mu(M)$ 
of volume-preserving diffeomorphisms of $M$. 
Arnold showed that fluid motions (in analogy with the classical case of the rigid body) 
in fact correspond to geodesics of the right-invariant metric on $\mathscr{D}_\mu(M)$ 
defined by the kinetic energy. 
This is a direct consequence of the least action principle and the postulate that fluid particles 
are allowed neither to fuse nor to split. 
Indeed, assuming appropriate smoothness conditions, 
let  $\gamma = \gamma(t,x)$ be the flow of the velocity field $v$, i.e. 
$$ 
\frac{d}{dt} \gamma(t,x) = v(t,\gamma(t,x)), 
\quad 
\gamma(0,x)=x. 
$$ 
Differentiating both sides of the flow equation in $t$ and using \eqref{ideal} 
leads immediately to the following second order system 
\begin{equation} \label{eq:geodesics} 
\frac{d^2 \gamma}{dt^2} 
= - \nabla p \circ \gamma \,,
\end{equation} 
which, roughly speaking, expresses the fact that the acceleration of the fluid is $L^2$-orthogonal 
(in the kinetic energy metric) to the space of divergence-free velocity fields. 
The latter constitute the tangent space at the identity to the group $\mathscr{D}_\mu(M)$ 
and 
the orthogonality condition represents the fact that the particle trajectories 
describe a geodesic curve in $\mathscr{D}_\mu(M)$.

More generally, let the fluid domain be an $n$-dimensional Riemannian manifold $M$ 
and let $\mu$ be the Riemannian volume form. The kinetic energy metric on $\mathscr{D}_\mu(M)$ 
is given at the identity $e$ by the $L^2$ inner product 
\begin{equation} \label{eq:L2metric} 
\langle v, w \rangle_{L^2} 
= 
\int_M \big( v(x), w(x) \big)_{T_xM} \,\mu 
\end{equation} 
of divergence-free vector fields $v, w \in T_e\mathscr{D}_\mu(M)$. 
As before, the trajectories of fluid particles satisfy the equations \eqref{eq:geodesics} 
and their velocities satisfy the Euler equations \eqref{ideal} 
with the nonlinear term $v{\cdot}\nabla v$ replaced now by the covariant derivative 
$\nabla_{\displaystyle v}v$ on $M$ 
and $\nabla p$ by the corresponding Riemannian gradient on $M$, see \cite{ak}. 
(This construction also applies if $\mu$ is an arbitrary volume form 
which does not coincide with the Riemannian volume form on $M$, 
provided that the divergence $\mathrm{div}\, v$ is taken with respect to $\mu$.)

\section{Ramifications of the Euler equations} 


\subsection{The Euler equations  with sources and sinks}
\label{s:geodesicflow} 
Various interesting, physically relevant and as yet unresolved, problems can be formulated already at this stage.
\begin{problem} 
Find a variational (preferably geodesic) formulation describing the motion of an ideal fluid 
in a fixed domain $M$ containing sources and sinks. 
What is the correct formulation of the variational problem: 
should one take into account the exterior forces and/or the ``memory" of the fluid?
\end{problem} 

For example, consider the case of a horizontal pipe with a fluid entering at one end and exiting at the other end. 
Such problems have a long history. 
On the one hand, as any mechanical system, fluid motions should obey some least action principle, 
see \cite{MatMet}. 
On the other hand, the energy of the fluid in the pipe may not be conserved since, 
depending on the boundary conditions, it could be supplied or drained at the two ends.
For instance, 
in addition to the equations \eqref{ideal} and the initial condition $v(0)$, 
the full system of the Euler-type equations in the 2D setting 
would have to include as data two other items: 
the function $v\cdot n$ describing the normal component of the velocity $v$ 
on the penetrable boundary, 
as well as 
the vorticity function $\omega:={\rm curl }\, v$ defined on the source part of the boundary 
through which the fluid is supplied. 
(Since vorticity is transported by the flow, this data will be sufficient
to define it for all times, see \cite{gie, yu, weigant, sueur}.) 

V.~Yudovich used such data to formulate a stability criterion for a steady pipe flow in 2D, 
which may hint at the appropriate boundary conditions needed to obtain a variational formulation, 
see \cite{yu, yu95}. 
We should add that for a ``dual" problem involving a fixed amount of fluid 
in a domain with a dynamic boundary its  Hamiltonian formulation is described in \cite{Marsden}. 


\subsection{The Euler equation for multiphase fluids and groupoids}
\label{s:groupoids} 

While Arnold's approach to fluids is limited to systems whose symmetries form a Lie group, 
there are many problems in fluid dynamics, such as free boundary problems, a rigid body in a fluid 
or fluid flows with vortex sheets, whose symmetries should instead be regarded as a {\it Lie groupoid}.
Groupoids can be thought of as groups with partially defined multiplication: for instance, fluid configurations with free boundary correspond to diffeomorphisms from one fluid domain to another; only maps for which the image of one coincides with the source of another admit composition (``multiplication").

In \cite{IK18, IK22} Arnold's framework was extended from Lie groups to  Lie groupoids to give a groupoid-theoretic description 
for incompressible  multiphase fluids, generalized flows, and fluid  flows with vortex sheets (the latter are flows whose  velocity field has a jump discontinuity along a hypersurface). 
A multiphase fluid consists of several fractions that can freely penetrate through each other without resistance and are constrained only by the conservation of total density. Beyond the vortex sheet setting, multiphase fluids arise e.g. in plasma physics and chemistry. Of particular interest are multiphase fluids with continuum of phases (or generalized flows), introduced by Brenier \cite{Brenier99}. One can think of them as flows in which every fluid particle spreads into a cloud thus moving to any other point of the manifold with certain probability \cite{sh}.

The Euler equations for multiphase flows on a Riemannian  manifold $M$ have the form
$$
\begin{cases}
\partial_t v_j + v_j\cdot \nabla v_j = -\nabla p\,, \\ 
  \partial_t \rho_j + {\rm div}\,(\rho_j v_j) = 0\,.
\end{cases}
$$
Here $\rho_1, \dots, \rho_n\in C^\infty(M)$ are  densities of $n$ phases of the fluid subject to the total incompressibility condition $\sum_{j=1}^n \rho_j = 1$, the vector fields $v_1, \dots, v_n\in \Vect(M)$ are the corresponding fluid velocities, and the pressure $p \in C^\infty(M)$ is common for all phases. For generalized flows the integer index $j=1,\dots,n$  enumerating the phases is replaced by a continuous parameter. In the case of vortex sheets the  densities are indicator functions of different parts of the manifold.

It turned out that in all of the above cases the corresponding configuration space has a natural groupoid structure. 
Using the corresponding Lie groupoids of multiphase  diffeomorphisms instead of the Lie group of 
volume-preserving transformations in Arnold's setting one  can describe the corresponding Lie algebroids and obtain geometric and Hamiltonian interpretations for the motion of the corresponding multiphase fluids, ``homogenized" vortex sheets, and 
generalized flows. Solutions of the above Euler equations were proved to be precisely the geodesics of an $L^2$-type right-invariant (source-wise) metric on the corresponding Lie groupoids of multiphase volume-preserving diffeomorphisms \cite{IK18, IK22}. Another interesting domain for applications of Lie groupoids is provided by elasticity theory, cf. 
\cite{mar-hugh, Kupf, JimLeon}.

\begin{problem}
Extend the geodesic and Hamiltonian descriptions of the Euler-Arnold equations on Lie groupoids  to problems of elasticity theory.
\end{problem} 

Many other open problems discussed below for the Euler equations related to diffeomorphism Lie groups can  also be posed for  the corresponding Lie groupoids, see e.g. Sections \ref{sec:cauchy}, \ref{sec:exp}, and \ref{sect:vortex}. For instance, it is natural to extends Arnold's study of the differential geometry of infinite-dimensional groups to those groupoids in view of possible applications to fluid stability problems.
\begin{problem} 
Describe the differential geometry (including computations of sectional and Ricci curvatures, conjugate points, etc.) 
for the right-invariant (source-wise) $L^2$-metric on the Lie groupoid of multiphase or generalized volume-preserving diffeomorphisms (analogous to Arnold's description of the differential geometry of the group $\mathscr{D}_\mu(M)$). 
\end{problem} 

\section{Variational setting for compressible fluids}
\subsection{Variational setting for shocks}\label{sect:shocks}
The inviscid Burgers equation 
$$
\partial_t v + v \cdot \nabla v=0
$$
describes freely moving non-interacting particles in a manifold $M$ of any dimension.
It  can  be also viewed as a geodesic equation, which in this case is defined 
on the full diffeomorphism group $\Diff(M)$ equipped with a non-invariant $L^2$-metric \cite{otto}. 
Once shock waves develop, the Lagrangian representation breaks down in the sense that 
the equation ceases to define an evolution in $\Diff(M)$; see e.g. \cite{KM07}. 
However, particles that stick inside the shocks continue to move along their own trajectories. 
For potential solutions with convex potentials there is a pointwise variational principle 
described by a ``circle law" proposed by Bogaevsky, see \cite{Bo} and its generalization in \cite{KSb, KSb12}. 
It prescribes the velocity $v^*$ of the common point of several colliding waves with velocities $v_i$: 
the joint velocity $v^*$ of the shock is given by the center of the smallest ball (a disk in 2D) 
covering all the velocities $v_i$ of the colliding waves. 

It is an interesting problem to formulate a more general variational principle for maps of $M$ to itself
(one should possibly consider Lipschitz maps to ensure differentiability almost everywhere) 
describing trajectories of particles, which would be valid before and after the formation of shocks 
and which would agree with both the non-invariant $L^2$-metric for the Burgers equation (before the collision) 
and the ``circle law" for particles sticking to each other after the collision.

\begin{problem} 
Is it possible to extend Arnold's geodesic framework 
from the group of diffeomorphisms $\Diff(M)$ to the semigroup of maps ${\mathcal Map}(M)$ 
that would capture both smooth solutions and their continuations beyond emergence of shock waves 
for Burgers-like and compressible fluid equations? 
\end{problem} 

To describe the motion of particles which are fused together inside shocks
one might employ the setting of generalized solutions (see \cite{ak, Brenier99, Brenier13})  to the Euler equations and the methods of control theory, 
which are well-adapted to study non-uniqueness of trajectories of dynamical systems. 
\bigskip

\subsection{Variational setting for sticking particles} 
There are various promising approaches to the variational formulation of the problem for sticking particles, see e.g.  
\cite{Brenier13, Sinai}. Here, we propose to look at it from yet another point of view.
We begin with the simplest situation: 
a motion of two sticking particles of equal mass moving without friction along a line. 
After the collision they form a new compound particle whose total mass and momentum are conserved. 
The kinetic energy obviously decreases upon collision. 
This loss can be interpreted as a transfer of energy to new unobservable degrees of freedom. 
For example, we can imagine that the particles move along two very close parallel lines 
and that, at the moment of their near-collision, they are joined by a rigid rod. 
The compound particle (in the shape of a dumbbell) will remain in the state of rotation: 
the angular coordinate of the axis of the rotating dumbbell is the new degree of freedom. 
Thus, a portion of the apparently vanishing energy has been allocated to this ``invisible'' degree of freedom. 
There may be physically different realizations of such invisible degrees of freedom. 
However, we only need to know that they exist and we are free to use them at will. 

Let $x_1, x_2$ be the coordinates of the particles with $x_1\le x_2$. 
The configuration space of our system is the half-plane $X=\{(x_1, x_2)|\  x_1 \le x_2\}\subset \R^2$. 
Let $\Delta=\{(x_1,x_2)|\  x_1=x_2\}$ be the diagonal and let $\Delta^\perp =\{(x_1,x_2)|\  x_1+x_2=0\}$ 
be its orthogonal complement. 
\smallskip

Extend the configuration space to the set $Z\subset\R_1^2\oplus\R_2^2$, where
\begin{align*} 
Z
&= 
(X\oplus\{0\})\cup ( \Delta\oplus\Delta^{\perp})
\\ 
&= 
\big\{ (x_1,x_2,y_1,y_2)|\ x_1\le x_2, y_1=y_2=0 \big\} 
\cup 
\big\{ (x_1,x_2,y_1,y_2)|\ x_1=x_2, y_1+y_2=0 \big\} 
\end{align*} 
is the union of the original space $X\subset \R^2_1$ and the plane $\Delta\oplus\Delta^\perp$ 
and where $\Delta\oplus\{0\}$ is identified with $\Delta\subset\R^2_1$. 

Now, consider two points $z_0\in X\subset Z$ and $z_1\in Z$ 
and a trajectory $z(t)$ in $Z$ for $0\le t\le 1$ with $z(0)=z_0$, $z(1)=z_1$ 
such that the action $J(z(\cdot))=\int_0^1\fr{1}{2}|\dot z(t)|^2 dt$ is minimal among all trajectories in $Z$ 
connecting $z_0$ and $z_1$. 
Let $P$ be the projection of $\R^2\oplus\R^2$ onto the first summand and define the trajectory $x(t)=P z(t)$. 
It is easy to see that $x(t)=(x_1(t), x_2(t))$ represents the motion of two particles on the line 
colliding and sticking upon collision with the total momentum being constant. 
Thus, we have established the variational principle in the simplest case of two particles on the line.

In a similar way we may consider a configuration of $n$ particles $x_1, \ldots, x_n$ on the line 
where $x_1\le x_2\le \ldots \le x_n$. 
Let $X\subset\R^n$ be the set of such configurations.  The set 
$X$ is stratified: 
let $\Delta_{m_1, \ldots, m_k}$ denote the set of $(x_1, \ldots, x_n)$ such that 
$$ 
x_1 = \dots = x_{m_1}, 
\  
x_{m_1+1} = \dots = x_{m_1+m_2}, 
\ 
\ldots 
\ 
x_{m_1 + \dots + m_{k-1}+1} = \dots = x_n 
$$ 
where $m_1+\ldots+m_k=n$. 
In $\R^n\oplus\R^n$ let 
\begin{align*} 
Y_{m_1, \ldots,m_k} 
&= 
\Delta_{m_1,\ldots, m_k}\oplus\Delta_{m_1, \ldots, m_k}^\perp 
\\ 
&= 
\bigg\{ 
(x_1, \dots, x_n, y_1, \ldots, y_n)\ |\ (x_1, \ldots, x_n)\in \Delta_{m_1,\ldots,m_k}, 
\begin{split} 
&y_1 + \dots + y_{m_1} = 0, \ldots 
\\ 
&y_{m_1+\dots+m_{k-1}+1}+\dots+y_n=0 
\end{split} 
\bigg\} 
\end{align*} 
and define the extended configuration space 
$$ 
Z 
= 
\big( X\oplus\{0\} \big) 
\bigcup 
\left( \bigcup_{\substack{k=1,\ldots, n-1 \\ m_1+\ldots+m_k=n}} Y_{m_1,\ldots,m_k} \right). 
$$ 
Let $P$ be the orthogonal projection from $\R^n_1\oplus\R^n_2$. 
Let $x_0\in X, \ z_0=x_0\oplus\{0\}$ and $z_1\in Z$. 
We can define the trajectory $z(t)\in Z, \ 0\le t\le 1$, connecting $z_0$ and $z_1$ 
and whose least action $J(z(\cdot))$ is minimal. 
Then, the trajectory $x(t)=Pz(t)$ in $X$ connects $x_0$ and $x_1=P z_1$ and describes the motion of 
sticking particles with the momentum preserved upon every collision. 

\bigskip

Now consider a continuum of material points distributed on the line. To be specific, consider the following situation: 
let $S$ be the segment $0\le s\le 1$ on the $s$-axis where $s$ is the label of a fluid particle. 
A fluid configuration is defined by the coordinate $f(s)$ for a particle with the label $s$, 
i.e. it  is a map $f: \ S\to\R$, 
and we assume that $f$ is a monotone function, 
i.e. $s_1\le s_2\Rightarrow f(s_1)\le f(s_2)$. 
The configuration space is $X=\{f\in W:=L^2(S,\R)\ | \ f \text{ is a monotone function on } S\}$. 
This space is stratified in the following way. 
Let $f(s)\in X$. This function may be constant 
on at most countably many intervals $\sigma$. 
Let $\Sigma$ be the collection of such intervals where $f(s)=\text{const}$ 
and 
define the stratum $X_\Sigma$ to be the set of all functions $f$ in $X$ such that 
$f|_\sigma=\text{const}$ for every $\sigma\in\Sigma$. 
Let $W_\Sigma=\{f\in W \ |\ f=\text{const} \text{  on each  }\sigma \in\Sigma\}$; 
then $X_\Sigma=X \cup W_\Sigma$. 
We see that its orthogonal complement is
$$
 W_\Sigma^\perp 
 = 
 \Big\{ 
 f\in W \ |\ \int_\sigma f \, ds = 0 \text{  for every  }\sigma\in\Sigma 
 \,\text{  and  }\, 
 f=0 \text{  outside}\, \cup_{\sigma\in\Sigma}\sigma 
 \Big\}. 
$$

Consider the space $W\oplus W$ and define the set $Z\subset W\oplus W$ 
(which is the desired extension of the space $X$) 
as follows.  
First, let $Y_0=X\oplus \{0\}\subset W\oplus W$. 
Next, for every $\Sigma\ne\emptyset$ let $Y_\Sigma=X_\Sigma\oplus W_\Sigma^\perp$ 
and set 
$$
Z := Y_0\, \bigcup\, \left( \bigcup_\Sigma Y_\Sigma \right) . 
$$ 

We are now in the position to formulate the variational problem whose solutions describe 
motions of a continuous family of particles on the line that stick upon collisions. 
Let $f_0 \in X$ be the initial position of the particles and set $g_0=f_0\oplus\{0\}\in Y_0\subset Z$. 
Let $g_1\in Z$. 
For any trajectory $g_t \in Z$ with $0\le t\le 1$ we define its action $J\{g_t\}=\int_0^1\frac{1}{2} |\dot g_t|^2 dt$. 
Among all trajectories $g_t$ in $Z$ connecting $g_0$ and $g_1$ we choose the one with minimal $J\{g_t\}$ 
(such trajectory exists and is unique). 
Now, define the trajectory $f_t$ in $X$ by $f_t=P\,g_t$, 
where $P$ is the projection from $W\oplus W$ onto the first coordinate. 
This is the desired trajectory of the system. 

\begin{remark} 
{\rm 
Note that the additional variables (the second term in $W\oplus W$ in the definition of the extended space $Z$) 
are necessary to define a sufficiently wide class of motions. If we simply defined $f_0, f_1 \in X$ 
then the action minimizing trajectory $f_t$ from $f_0$ to $f_1$ would be a motion along the straight segment 
connecting these two points (since $X$ is a convex set). For such a trajectory the particles would not collide at all 
(or, if you wish, they would collide only at $t=1$). 
} 
\end{remark} 

\begin{remark} 
{\rm 
If we attempt to use this approach to construct the motion of a continuum of particles in $\R^d, \ \ d>1$ 
then we encounter a new difficulty: the set of admissible configurations of the particles is no longer convex. 
Therefore, particle collisions are not as easily controlled and parametrized as in the one-dimensional case. 
In particular, it is unclear how to use this method to give a variational description of the formation of ``shock waves'', 
i.e., hypersurfaces where the mass is concentrated with positive density. 
In the theory of shock waves individual particle motions are described for potential solutions only, 
cf. Section \ref{sect:shocks}. 
 
However, there exists another class of sticking flows in $\R^3$, namely, flows with constant density 
(and decreasing energy), see \cite{sh2}. 
Such flows are dissipative weak solutions of the Euler equations and are vaguely similar to turbulent flows. 
Their construction is based on different ideas (not on a variational principle) and it would be interesting to define them 
in a way corresponding to the above 1D systems. 
} 
\end{remark} 
We thus arrive at the following problem: 
\begin{problem} 
Find a variational description of the system of material particles, moving in $\R^d$ with $ d>1$ and sticking upon collision, 
for both potential and non-potential fluid flows. 
This formulation should be sufficiently flexible to describe the formation and development of shock waves 
in the system described by the multi-dimensional Burgers equation. 
\end{problem} 
\begin{remark} 
{\rm 
The above variational principle should be closely related to the models of adhesion
particle dynamics studied in \cite{Sinai}. It seems to give the same results  for a finite number of sticking particles. 
However, the approach in \cite{Sinai} is fundamentally one-dimensional, while the approach described above 
may be extended to higher dimensions as well. 

Indeed,  for a finite number of particles even in higher dimensions one can assign which particles stick together. 
However, on the way to this collision they may bump into other particles. 
For a finite number of particles there is an excuse that the probability of this happening is zero. 
But, in the case of a continuum, e.g., given a continuous density at the initial moment particle trajectories will intersect en masse. 
One might allow that by envisioning a ``dusty matter" with a multi-flow structure and no pressure 
(for example, this might be the case of stars in the universe when different streams of stars move in different directions 
in the same volume). 
Alternatively, one might confine to piecewise-smooth flows with stratified density 
(supported on a stratified manifold with components of different dimensions). 
The latter setting is close to Kantorovich's theory of optimal mass transport, cf. \cite{Brenier99, Brenier13}, 
where such a motion in 1D with gluing of particles and shock waves is described 
and the variational principle is written in terms of differential inclusions on the space of transport maps. 
} 
\end{remark} 

\bigskip 

\section{Rigid and fluid modelling of invariant metrics} 
Since the work of V.~Arnold it is well-known that 
while the Euler equation of a rigid body corresponds to a {\it left-invariant metric} 
(depending on the body shape) on the group $SO(3)$, 
the Euler equations of ideal hydrodynamics correspond to an $L^2$ {\it right-invariant  metric} 
on the diffeomorphism group $\mathscr{D}_\mu(M)$, see \cite{arn66, MatMet}. 
Left invariance of the rigid body metric is related to the fact that the body's energy 
depends on the angular momentum in the body and does not depend on its position in the ambient space. 
On the other hand, right invariance of the fluid metric is related to the norm of the velocity field in the space 
but does not depend on the parametrization of fluid particles. In other words, the energy metric on $SO(3)$ is left-invariant because the space $\R^3$ is {\it isotropic}, while the metric on $\mathscr{D}_\mu(M)$ is right-invariant because the fluid is 
\emph{homogeneous}. 
We thus see that the reasons of right- and left-invariance are quite different.

A natural question arises: 
does there exist an interesting mechanical system with $SO(3)$ as a configuration space 
such that the energy metric is {\it right}-invariant? The answer is {\it yes} and here is its description.

First, we define an object called a ``hedgehog". This is a ball $B$  whose center is a fixed point $O\in \R^3$ 
so that it can freely rotate around it. Suppose that at every point on the surface of the ball grows a ``needle", 
i.e., there is a (sufficiently long) radial segment. The whole structure rotates around the fixed center $O$ 
as a solid body so that its configuration space is $SO(3)$.

Next, suppose that on every needle there is an infinitesimal point mass which is able to move freely along the needle; 
let us call it a bead. Let $\rho(\omega)$ be the angular mass density, so that the bead mass in the solid angle $d\omega$ 
is $\rho(\omega) d\omega$ (the hedgehog itself is massless). 
Now, suppose that there is a closed surface $S$ surrounding the ball so that every needle pierces $S$. 
Lastly, suppose that every bead is forced to remain at all times on the surface $S$ 
and  
at the same time is confined to its own needle. 
It is natural to call this system ``beads on the hedgehog" (we do not consider here its practical realizations). 
This system depends on the angular bead density $\rho(\omega)$ and the surface $S$.

\begin{figure}
\centering
\includegraphics[width=2.5in]{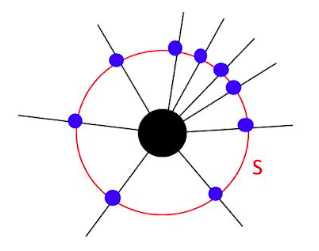}
\vspace{-7pt}
\caption{\small  Beads on a spherical ``hedgehog model": the angular density of beads on a concentric sphere $S$ is arbitrary. The corresponding metric on $SO(3)$ is left-invariant.}
\vspace{-10pt}
\label{Fig.0a}
\end{figure}

\begin{figure}
\centering
\includegraphics[width=2.5in]{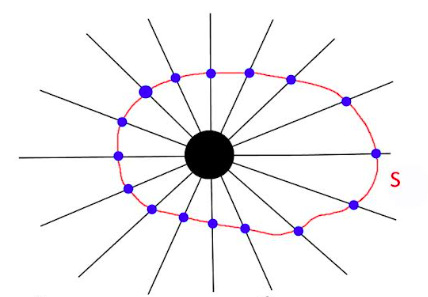}
\vspace{-7pt}
\caption{\small  Beads on a star-shaped ``hedgehog model": the surface $S\subset \R^3$ is arbitrary, the angular density of beads $\rho(\omega)$ is constant. The corresponding metric on $SO(3)$ is right-invariant.}
\vspace{-10pt}
\label{Fig.0b}
\end{figure}

\medskip

There are two cases of the hedgehog having some symmetry. In the first case, the surface $S$ is a sphere concentric with the hedgehog, while the bead density $\rho(\omega)$ is arbitrary (Figure~\ref{Fig.0a}). It is easy to see that this is the same as a solid body with a fixed point, the Euler top. In this case the metric on $SO(3)$ is left-invariant: the energy depends on the angular velocity of the system in the body but not in the space.

In the second case, $S$ is an arbitrary star-shaped surface (fixed in $\R^3$), while the angular density $\rho(\omega)$ 
is constant (Figure~\ref{Fig.0b}). In this case the metric on the group $SO(3)$ is right-invariant: 
now the energy depends on the angular velocity of the system in the space but not in the body. 
Note that the equations describing this dynamics will be the Euler-Arnold equations of the right-invariant metric on $SO(3)$ 
which differ only by a sign from the standard equations for the Euler top.

In all other cases the metric is neither left- nor right-invariant. It would be interesting to investigate this system 
for generic $S$ and $\rho(\omega)$. It would also be of interest to define and study a model system 
where the configuration space is $SL(2)$, since this group looks more ``liquid-like". 
A similar question about possible fluid models is very intriguing.
\bigskip

\begin{problem} 
Describe a model for an $L^2$ {\rm left-invariant metric} on the group $\mathscr{D}_\mu(M)$. 
Are there any interesting physical systems demonstrating that type of invariance? 
\end{problem}

\section{The Cauchy and boundary value problems} \label{sec:cauchy}
%
\subsection{On local well-posedness of the Cauchy problem} 
First rigorous results on local existence and uniqueness of solutions to the Cauchy problem 
for the incompressible Euler equations \eqref{ideal} were obtained in the 1920s 
by Gunther \cite{gu} and Lichtenstein \cite{li} in the class of H\"older $C^{1,\alpha}$ spaces. 
Global existence in 2D was established shortly thereafter by Wolibner \cite{wo}. 
Various subsequent extensions and improvements of these results in H\"older, Sobolev $H^s$ and $W^{s,p}$ 
and more exotic Besov $B^s_{p,q}$ and Triebel-Lizorkin $F^s_{p,q}$ spaces 
can be found in the papers \cite{yu, ka, ema, bkm, kap}; 
see also recent monographs and surveys \cite{mb, bcd, co, bt0, de}. 

Roughly speaking, a Cauchy problem is \textit{locally well-posed} in a Banach space $X$ 
(in the sense of Hadamard) 
if given any initial data in $X$ there exists a $T > 0$ and a unique solution in the space $C([0, T], X)$ 
which depends at least continuously on the data. 
Otherwise, the problem is said to be \textit{ill-posed} in $X$. 
A number of ill-posedness mechanisms have been investigated in the literature, 
from loss of regularity properties of the solution map, to energy decay, to nonuniqueness 
and finite time blowup. 
Although global (in time) well-posedness of the 3D Euler equations has long been seen as the major open problem 
in analysis and PDE, 
interesting questions concerning local well-posedness (in any dimension) 
have also remained open for a long time. 

Recall that the Cauchy problem for the incompressible Euler equations is not well-posed 
in the standard H\"older spaces in the sense that solutions may depend discontinuously 
on general initial data in $C^{1,\alpha}$. However, this dependence is known to be continuous 
in, e.g., the ``little" H\"older space (essentially, the completion of smooth functions in the H\"older norm), 
as well as in $W^{s,p}$ Sobolev spaces with $p \geq 2$ and $s>n/p+2$, 
cf. e.g., \cite{ema, kap0, my1}. 
Somewhat more refined existence and uniqueness results 
are available in $B^1_{\infty,1}$ and $B^{n/p +1}_{p,1}$ where $1<p<\infty$, 
see e.g. \cite{pp, vi, chae}. 
%
On the other hand, there are examples of 3D solutions of \eqref{ideal} which exhibit instantaneous loss of regularity. 
Examples in $C^\alpha$ with $0<\alpha<1$, as well as in $B^1_{\infty,\infty}$, $F^1_{\infty, 2}$ 
and $\log\mathrm{Lip}^\alpha$ with $0<\alpha \leq 1$ 
were constructed in \cite{bt, l-r,  my}. 
Many other ill-posedness results 
in the borderline Sobolev and Besov spaces 
$W^{n/p+1,p}$, $B^{n/p +1}_{p,q}$ with $1 \leq p < \infty$ and $1 < q \leq \infty$ 
and in the classical spaces $C^k$ and $C^{k-1,1}$ can be found in \cite{bourg-d, bourg-d-2, elmas}. 
%

One remaining case of particular interest can be formulated as follows. 
\begin{problem} 
Are the Euler equations \eqref{ideal} ill-posed in the Besov spaces $B^1_{\infty, q}$ for $1<q<\infty$? 
\end{problem} 
%
To put this problem in a functional space context, first recall that 
$$ 
H^s 
\subset\, 
C^{1+\alpha} 
\subset\, 
B^1_{\infty, 1} 
\subset\, 
C^1 
\subset\, 
\mathrm{Lip} \, 
\subset\, 
F^1_{\infty, 2} 
\subset\, 
B^1_{\infty,\infty} 
\subset\, 
\mathrm{logLip} \, 
\subset\, 
C^\beta 
$$ 
for any $0< \alpha, \beta <1$  and $s>n/2 + 1+\alpha$ 
and 
next observe that, in the specified range, the $B^1_{\infty,q}$ spaces interpolate between 
$B^1_{\infty,1}$ and the Zygmund space $B^1_{\infty,\infty}$. 
%

%
%
%

\subsection{A two-point boundary value problem on the diffeomorphism groups} 
Turning to the geodesic equation \eqref{eq:geodesics}, 
consider the Sobolev completion $\mathscr{D}^s_\mu(M)$ of the diffeomorphism group $\mathscr{D}_\mu(M)$. 
As is well known, the Cauchy problem for \eqref{eq:geodesics} can be solved (for small values of $t$) 
by standard Banach contraction arguments provided that $s>n/2+1$, 
cf. \cite{ema}. 
Consequently, the $L^2$ metric \eqref{eq:L2metric} admits a smooth Riemannian exponential map 
\begin{equation} \label{eq:exp} 
\exp_e: T_e\mathscr{D}_\mu^s(M) \to \mathscr{D}_\mu^s(M) 
\end{equation} 
defined in a neighbourhood of the zero vector by 
$\exp_e tv_0 = \gamma(t)$ 
where $\gamma$ is the unique geodesic from the identity element $e$ with initial velocity $v_0$. 
Furthermore, as in the classical finite-dimensional Riemannian geometry, 
we have $d\exp_e(0)=\mathrm{id}$, 
and therefore $\exp_e$ is a local diffeomorphism of Banach spaces by the inverse function theorem. 
In particular, this implies   local well-posedness
of the Euler equations in $H^s$ for $s>n/2+1$, 
as well as 
unique solvability of the two-point boundary value problem for the geodesic equation \eqref{eq:geodesics} 
in any sufficiently small neighborhood of $e$ in $\mathscr{D}^s_\mu(M)$. (In geometric language 
Wolibner's global existence and uniqueness result in \cite{wo} amounts to a statement about 
geodesic completeness of the manifold $\mathscr{D}_\mu^{1,\alpha}(M)$ of H\"older diffeomorphisms 
under the right-invariant $L^2$ metric \eqref{eq:L2metric} when $n=2$.)

A natural question is whether the two-point problem holds in the large. 
We formulate it as two related problems. 
\begin{problem}[Two-point boundary value problem] \label{2ptbvp} 
Let $M$ be a compact two-dimensional Riemannian manifold 
and let $\varphi$ be a volume-preserving diffeomorphism of $M$ of Sobolev class $H^s$. 
\begin{enumerate} 
\item[(i)]{\rm (Surjectivity Problem)} 
Find a divergence-free vector field $v \in H^s(M)$ such that $\exp_e{v} = \varphi$. 
\item[(ii)]{\rm (Variational Problem)} 
Find a curve $\gamma(t)$ in $\mathscr{D}_\mu^s(M)$ from $e$ to $\varphi$ 
which minimizes the $L^2$ energy functional 
$ 
\mathcal{E}(\gamma) = \frac{1}{2} \int_0^1 \| \dot{\gamma}(t)\|^2_{L^2} dt. 
$ 
\end{enumerate} 
\end{problem} 
Although two-point boundary value problems in hydrodynamics are no less fundamental than 
the Cauchy problem they have not received as much attention. 
From the geometric point of view $(i)$ and $(ii)$ may be regarded as infinite-dimensional versions 
of the classical Hopf--Rinow theorem. 
%
%
One strategy for $(i)$ is to follow the classical argument of 
Hopf and Rinow compensating for the lack of local compactness with {\it a priori} estimates derived 
with the help of weak solutions and Lyapunov functions \cite{sh1}. 
Another approach could use the properties of the exponential map as a nonlinear Fredholm and quasiruled map 
\cite{emp, shn0, sh3}. 
In connection with $(ii)$ we mention a surprising result of  \cite{shn, sh} 
that there exist volume-preserving diffeomorphisms of a simply connected compact 3-manifold 
which cannot be joined by a shortest path in $\Diff_\mu^s(M)$. 
For the two dimensional case partial results can be found in \cite{mp}. 

\subsection{Global geometry of the group of volume-preserving diffeomorphisms}

Now we turn to questions concerning the global geometry of the group $\mathscr{D}_\mu^s$. 
The following problem is related to Problem \ref{2ptbvp}. 
\begin{problem} 
Does the energy (action) functional corresponding to the $L^2$ kinetic energy metric \eqref{eq:L2metric} 
in 2D hydrodynamics satisfy the Palais--Smale condition? 
\end{problem} 

Another interesting question is related to closed geodesics. 
Suppose that $M$ is a compact surface (possibly with boundary) of genus at least $2$ 
or that $M$ is a multi-connected bounded domain in $\mathbb{R}^2$ with at least two holes. 
\begin{problem} \label{prob-closed} 
Does there exist a closed geodesic in $\mathscr{D}^s_\mu(M)$? 
\end{problem} 
Note that the Kelvin--Helmholtz theorem implies that any geodesic loop in $\mathscr{D}^s_\mu(M)$ 
is necessarily a closed geodesic. 
One may try to construct a suitable Lyapunov function (see below) to show that in this case 
a geodesic never returns to its initial configuration. 

\medskip

Yet another basic problem concerns the fluid configuration space $\mathscr{D}_\mu(M)$ itself. 
As an infinite dimensional (Frechet) Lie group it can be viewed as a Riemannian homogeneous space 
equipped with a right-invariant $L^2$ metric \eqref{eq:L2metric}. 
However, the group $\mathscr{D}_\mu(M)$ is not a Riemannian symmetric space. This means that the (geodesic) central symmetry about the identity in $\mathscr{D}_\mu(M)$ is not an isometry of the $L^2$-metric. In this respect, we have to keep in mind that the geodesic symmetry is not the same as the group-theoretic symmetry, i.e. the map $g\to g^{-1}$. The latter (group-theoretical) inversion maps a right-invariant metric on the group into a left-invariant one, hence, indeed, it is not an isometry. However, on the group of finitely differentiable diffeomorphisms like $\mathscr{D}^s_\mu(M)$ the group inversion is not even differentiable, while if we use the Holder $C^{k,\alpha}$-class diffeomorphisms, it is not even continuous.  Unlike this, the geodesic central symmetry is a smooth map in $\mathscr{D}_\mu(M)$ for any reasonable model space (like $H^s$ and $C^{k,\alpha}$), but the lack of its   isometry property is not so immediately seen. Hence, phrasing our question somewhat informally, we may ask

\begin{problem} 
How ``far" is $\mathscr{D}_\mu(M)$ from being an infinite-dimensional (locally) symmetric Riemannian space?
\end{problem} 
\smallskip

This question is related to the following long standing (although little known) paradox. 
Consider the parallel sinusoidal steady fluid flow given by the stream function $\psi = \cos( k\,y)$
on the two-dimensional torus. Then well-known Arnold's theorem claims that the sectional curvature of the group 
of exact area-preserving torus diffeomorphisms is nonpositive in all  (and negative in most) two-dimensional directions containing the direction given by $\psi$, see \cite{arn66, ak}. (There is a similar statement for 
a plane-parallel flow in a periodic channel.) 
 Following  Arnold's idea on an intrinsic relation between  negative curvature and the flow (Lagrangian) instability, one could expect that {\it any} plane-parallel flow is unstable. But this does not seem to be the case, since there are 
 (Eulerian) stable parallel flows (for instance those with convex velocity profiles on ``short tori", see discussion in \cite{ak}), while  Lagrangian and Eulerian instabilities are closely related, cf. \cite{pr04}.
 
The root of this misunderstanding lies in our ``symmetric" intuition. In fact, this relation between the curvature sign 
and (in)stability of geodesics exists for  symmetric spaces (say, on  groups with bi-invariant metrics), see e.g. \cite{Miln, ak}. On the other hand, the group $\mathscr{D}_\mu(M)$ is not symmetric, but rather ``chiral": as a Riemannian space it is somewhat twisted in one direction, and hence we observe a discrepancy between instability of geodesics and its negative curvature. The chirality of $\mathscr{D}_\mu(M)$ might have
some other, more profound, consequences beyond the instability issues, which would be interesting to explore.

\bigskip


\section{Partial analyticity of solutions in 2D}

\subsection{Analyticity of particle trajectories}\label{sect:analytic}
The Euler equations keep bringing surprises 
-- such as the following relatively recent theorem. 
Suppose that $M$ is a compact 2D real analytic manifold or a bounded domain with analytic boundary. 
Let $u$ be a solution of the Euler equations in $M$ of Sobolev class $H^s$ for $s>2$, 
obeying the slip condition $u\,\|\,\partial M$ if $\partial M \ne \emptyset$. 
Finally, let $x(t)$ be any particle trajectory satisfying the flow equation $\fr{d}{dt}x(t) = u(t,x(t))$. 

\begin{thm}
Any particle trajectory
$x(t)$ of the flow is a real-analytic function of $t$.
\end{thm}

This theorem was first established by Serfati \cite{Ser}, 
then independently by  Inci, Kappeler, and Topalov \cite{IKT2013}, 
Shnirelman \cite{shn-a}, 
Constantin, Vicol and Wu \cite{cvw}, 
Zheligovsky and Frisch \cite{zh-f} and, in case of stationary flows, Nadirashvili \cite{nadir}. 
The proofs by the above authors are based on two entirely different ideas. 
In the works \cite{Ser, cvw, zh-f}
the function $x(t)$ was formally expanded in the Taylor series, 
whose convergence was proved using commutator estimates. Thus, this was a straightforward proof. 
On the other hand, the works \cite{IKT2013, shn-a, nadir} 
used the Lagrangian description of the fluid motion, 
i.e. they considered the flow as the motion $g_t$ on the infinite-dimensional manifold $\mathscr{D}_\mu^s(M)$, 
equipped with the $L^2$-metric, along a geodesic. 
The manifold $\mathscr{D}_\mu^s(M)$ is a real-analytic Banach manifold and the geodesic spray 
(generating the geodesic flow in the tangent bundle to $\mathscr{D}_\mu^s(M)$) is an analytic vector field, 
since it can be extended to the complexification $\C \mathscr{D}_\mu^s(M)$ as a holomorphic vector field. 
Then the standard theorems of existence, uniqueness and analytical dependence of solution $g_t$ on  $t$ 
and on the initial condition $x_0$ can be applied, since they hold for any analytic Banach manifold \cite{hille-ph, dieu}. 

The work \cite{shn-a} was based on a similar idea: 
using the Kelvin--Helmholtz vorticity theorem, one can reduce the Lagrange equation to a first order equation 
of the form 
$\fr{d}{d t} g_t = V(g_t)$ on $\mathscr{D}_\mu^s(M)$, 
where $V$ is an analytic vector field on the infinite-dimensional manifold $\mathscr{D}_\mu^s(M)$. 
Then the above basic existence, uniqueness, etc. theorems are applicable, and the same result follows: 
the flow $g_t\in \mathscr{D}_\mu^s(M)$ is an analytic curve depending analytically on $g_0$. 

Interestingly, the latter work closely follows the original approach of L.~Lichtenstein \cite{li-0} 
in his proof of the local in time existence and uniqueness of solutions to the Euler equation. 
Lichtenstein proved that the vector field $V$ is Lipschitz (in fact, he proved that it is $C^1$) 
and that it can be continued analytically in the complexification of $\mathscr{D}_\mu^s(M)$. 
Thus, in 1925  he was just one step away from proving that the flow $g_t$ is analytic in $t$! 
However, Lichtenstein's work appeared  roughly 10 years before the Banach spaces acquired their name; 
about 20 years before the concepts of complex analysis (like the analytic implicit function theorem) 
were extended to the complex Banach spaces; 
and about 30 years before it was acknowledged that the basic concepts of smooth 
and analytic topology including the theory of ODEs can be transferred to the complex Banach manifolds. 
So, if he made that step his discovery would be truly extraordinary. 

\begin{problem}
Are particle trajectories analytic in time in any dimension? 
For fluids on manifolds with boundary how does this analyticity depends on whether the boundary is analytic or not?
\end{problem}

\medskip

Nadirashvili \cite{nadir} proved analyticity of flow lines of a {\it stationary} solution  to the 2D Euler equation. 
His theorem is local and holds independently of the analyticity (or the lack thereof) of the boundary $\partial M$. 
It follows the classical analyticity proof for solutions of analytic elliptic equations 
and uses the fact that an elliptic equation becomes hyperbolic if one of the variables, 
say $x_1$, is replaced by $\mathrm{i} x_1$. 
(It is worth recalling that flow lines and vorticity lines coincide for stationary 2D solutions 
and hence are analytic simultaneously.) 

\begin{remark}
{\rm
In dimensions $d=2 $ and $3$ the analyticity of particle trajectories in $\mathbb R^d$ (i.e. in the case without boundary) was proved in \cite{cvw}. Apparently the higher dimensional case does not present fundamentally new
difficulties, and the study in  \cite{cvw} was confined to  low dimensions because of their physical significance.
The paper \cite{Hern}  covers several related cases, including vortex patches among other situations.

There is also the following related problem: Are particle trajectories of Yudovich solutions analytic in time? (Recall that the Yudovich class consists of divergence-free vector fields $u$ satisfying $\|{\rm curl}\, u\|_{L^\infty}<\infty$. This inequality implies that the field $u$ has the Osgood property \cite{Petrovsky}, which guarantees the uniqueness of trajectories.)
It is currently only known that they are Gevrey regular, due to the result in \cite{Gam}, see also \cite{chem}.
It seems to be unknown if this Gevrey regularity is sharp.\footnote{We thank the anonymous referee for this remark.} 

The difficulty here  is as follows. It is known that Sobolev vector fields are integrated to Sobolev diffeomorphisms; in other words,  Sobolev vector fields form the Lie algebra of the group of  Sobolev diffeomorphisms. However, it is unknown what would be the result of integration for vector fields from the Yudovich class, i.e. what are ``Yudovich homeomorphisms". Do such homeomorphisms form a group? Does this ``group" admit the structure of a real Banach manifold? Can this ``manifold" be complexified? Presumably, the answers are ``no" to all those questions, and one has to look for other approaches to this problem.}
\end{remark}

In regard to analyticity, Lebeau \cite{Lebeau}
considered piecewise-continuous solutions of 2D Euler equations 
which are irrotational outside of a time-dependent curve and have a tangential discontinuity on the curve, 
cf. Section \ref{sect:vortex} on vortex sheets. 
Lebeau proved that the curve of discontinuity is analytic as long as such a solution exists, 
i.e. the vortex sheets are analytic in 2D.

\medskip
\subsection{Stationary flows and partially analytic functions}
For a stationary solution $u$ of the Euler equations particle trajectories are the same as flow lines. 
The stream function $\psi(x)$ of a stationary flow $u=\nabla^\perp\psi$ is a peculiar function: 
it may be an $H^s$ function but its level lines are analytic. 
A function $\psi(x)$ whose level lines are real-analytic will be called a {\it partially analytic function}. 
The set of partially analytic functions is by no means a vector space: 
just consider a pair of such functions $\psi_1$ and $\psi_2$ 
whose (analytic) level lines are transversal to each other; 
then the level lines of $\psi_1+\psi_2$ are not necessarily analytic. 
Thus, one arrives at an important problem to find a natural structure on the set of partially analytic functions.

One natural idea would be to choose  level lines, rather than the values at  points, 
as an adequate representation of a partially analytic function. 
For example, consider a function $\psi(x_1,x_2)$ defined on a periodic curvilinear strip 
$M=\{(x_1,x_2)\,|\, g(x_1)\le x_2\le h(x_1)\}$, where $g$ and $h$ are real periodic functions with the same period $\ell$, 
and such that every level line $\psi(x) = {\rm const}$ has an equation 
$x_2 = a(x_1, \psi)$ 
with $a(x_1,0) \equiv g(x_1)$ and  $a(x_1,1) \equiv h(x_1)$. One can assume that $x_1\in \mathbb{T} =\R/\ell\,\mathbb{Z}$.
Then the function $a(x_1,\psi)$ uniquely defines $\psi(x_1, x_2)$ in the flow domain $M$. 
The function $a(x_1,\psi)$ is of class $H^s$ and it is analytic in $x_1$ for any fixed $\psi$. 

To be more precise, let us define the corresponding function spaces. 
Let us fix $\sigma>0$.

\begin{defn}
The space $X^s_\sigma$ consists of real-analytic  functions $f(x_1)$, $x_1\in\mathbb{T}$ which can be analytically continued into the strip $|{\rm Im}\, x_1|\le \sigma$ and such that $f(\cdot\pm i\sigma)\in H^s(\mathbb{T})$, where the norm is $||f||_{X^s_\sigma}=\|f(\cdot - i\sigma)\|_{H^s(\mathbb{T})}+\|f(\cdot + i\sigma)\|_{H^s(\mathbb{T})}$.
\end{defn}

\begin{defn}
The space $Y^s_\sigma$ consists of functions $a(x_1,\psi)$ such that 

(i) for any $\psi\in[0,1]$, $a(x_1,\psi)\in X^s_\sigma$ ;

(ii) the functions $a(x_1\pm i\sigma,\psi)$ belong to $H^s(\mathbb T\times[0,1])$.

The norm in the space $Y^s_\sigma$ is defined as follows:
$$||a||_{Y^s_\sigma}=||a(\cdot+i\sigma,\cdot)||_{H^s(\mathbb{T}\times[0,1])}+||a(\cdot-i\sigma,\cdot)||_{H^s(\mathbb{T}\times[0,1])}$$
\end{defn}

And, at last, we define the space $Z^s_\sigma$:

\begin{defn}
A function $\psi(x_1, x_2)$ defined in the domain $M=\{(x_1,x_2)\,|\, g(x_1)\le x_2\le h(x_1)\}$ belongs to the space $Z^s_\sigma$ if its level lines $\psi={\rm const}$  can be defined by the equation $x_2=a(x_1,\psi)$, where the function $a(\cdot,\cdot)$ belongs to the space $Y^s_\sigma$. The norm in $Z^s_\sigma$ is induced from the space $Y^s_\sigma$.
\end{defn}

An immediate application of this space is to the description of the set of stream functions 
of stationary solutions of the Euler equations in the periodic domain 
$M=\{(x_1,x_2)\,|\, g(x_1)\le x_2\le h(x_1) \,\,{\rm for}\,\, x_1\in\mathbb{T} \}$. 
Indeed,  if we consider stationary solutions $u(x_1,x_2)$ in $H^s$ with stream functions $\psi(x_1,x_2)$ in $H^{s+1}$ 
then the set of stationary solutions does not form a smooth manifold in $H^s$. 
This difficulty was partially circumvented by Choffrut and Sverak \cite{chsv} 
by resorting to the $C^\infty$ Frechet spaces and the Nash-Moser-Hamilton inverse function theorem. 
In this $C^\infty$ 
setting the stationary solutions form a smooth manifold parametrized locally by distribution functions of the vorticity. 
However, those tools  might be too powerful  for the task in finite smoothness. 

On the other hand, using the space $Y^s_\sigma$, 
Danielski \cite{Danielski} established a local description of the set of stationary flows in the periodic channel. Let $u_0$ be the velocity field of a parallel flow satisfying several conditions in the domain $M_0=\mathbb{T}\times [0,1]$. Namely, assume that  $u_0(x_1, x_2)=(U_0(x_2),0)$ satisfies 
(1) $U_0(x_2)> 0$,  (2) $U_0'(x_2)=F_0(x_2)>0$, (3) $U_0''(x_2)>0$; and finally, 
(4) let $\psi_0(x_1, x_2)=\Psi_0(x_2)=\int_0^{x_2} U_0(t) d t$ be the stream function of the flow $u_0$, then $\psi_0$ satisfies the boundary conditions $\psi_0(x_1, 0)=0, \psi_0(x_1,1)=1$. Its level curves $\psi_0={\rm const}$ have the equation $x_2=a_0(\psi)$ where $a_0(\psi)$ is the inverse function to $\Psi_0(\cdot)$.

One can regard the above periodic domain $M=\{(x_1,x_2)\,|\, g(x_1)\le x_2\le h(x_1) \,{\rm for}\,\, x_1\in\mathbb{T} \}$ as 
being close to the parallel strip $M_0$. We are looking for the stationary flows in $M$ which are close to the parallel flow $u_0$ in $M_0$. 
\begin{thm} 
For any parallel flow $u_0$ possessing the properties {\rm (1)--(4)} there exists $\varepsilon>0$ such that the following holds.  Suppose that  $||g(x_1)-0||_{X^s_\sigma}<\varepsilon$ and $||h(x_1)-1||_{X^s_\sigma} <\varepsilon$. Then there exist stationary flows $u\in H^s$ close to $u_0$ in the following sense:
\begin{itemize} 
\item[(1)] 
for each solution $u$, its stream function $\psi\in Z^{s+1}_\sigma$;
\item[(2)] 
the stream function $\psi$ satisfies relation $\Delta\psi=F(\psi)$ for a certain monotone function $F\in H^{s-1}$ close to $F_0$ in $H^{s-1}$;
\item[(3)] 
the functions $\psi$ form an analytic submanifold $\Sigma\subset Z^s_\sigma$ which is locally analytically diffeomorphic to a neighbourhood of $F_0$ in $H^{s-1}[0,1]$.
\end{itemize} 
\end{thm}

The stationary flows with stagnation points present some additional difficulties.

\medskip

\subsection{An attractor of 2D Euler equations and its semianalytic structure}

Consider a compact analytic Riemannian 2-manifold $M$ with or without boundary, for example, a 2-torus. 
Let $YU(M)$ be the Yudovich space of vector fields on $M$ consisting of  divergence-free vector fields $u$ such that 
${\rm curl}\, u$ is in $L^\infty(M)$. 
The Euler equations define a perfect dynamics on the space $YU(M)$, 
i.e. a one-parameter group $S$ of transformations $S_t: YU(M)\to YU(M)$ continuous in the $H^1$ topology 
(i.e. weakly continuous in $YU(M)$). 
For any $u\in YU(M)$ let $\mathscr{O}(u)$ be the {\it orbit} of $u$. 
Let $\bar{\mathscr{O}}(u)$ be the weak closure of $\mathscr{O}(u)$ in $H^1(M)$. 
Note that for any $v\in \mathscr{O}(u)$, one has $||{\rm curl}\, v||_{L^2}=||{\rm curl}\, u||_{L^2}$ 
and for any $w\in\bar{\mathscr{O}}(u)$, $||{\rm curl}\, w||_{L^2}\le ||{\rm curl}\, u||_{L^2}$.
\begin{defn} 
A field $u\in YU(M)$ is called a {\rm generalized minimal flow} (or a GM-flow) 
if for any $w\in\bar{\mathscr{O}}(u)$ we have 
$\| {\rm curl}\, w \|_{L^2} = \| {\rm curl}\, u \|_{L^2}$. 
The set of all generalized minimal flows is denoted by $\mathscr{GM}$. 
\end{defn} 
This definition has the following meaning.  For any fluid flow its vorticity is transported by the fluid 
(the Kelvin--Helmholtz theorem). Thus, the vorticity field is deformed by the flow and, as $t\to\infty$, 
this deformation effectively leads to {\it mixing}. 
The mixing operator $K$ has the form $Kf(x)=\int_M K(x,y) f(y) \,dy$, 
where the kernel $K(x,y)$ is a non-negative measure in $M\times M$ such that $\int_M K(x,y) \,dx\equiv 1$ 
and $\int_M K(x,y) \,dy \equiv 1$ (i.e. $K$ is a bistochastic operator). 
Any mixing operator is a contraction in $L^2(M)$. 
Thus, for any $w \in \bar{\mathscr{O}}(u)$, ${\rm curl}\, w = K({\rm curl}\, u)$ for some mixing operator $K$. 
If $u\in{\mathscr{GM}}$ then for any $w\in \bar{\mathscr{O}}$, ${\rm curl}\, w$ is equimeasurable with ${\rm curl}\, u$. 
In other words, ${\rm curl}\, u$ is not mixed by the Euler flow; 
the measure $\mu=({\rm curl}\, u)\, dx$ can be disintegrated into components $\mu_\alpha$ 
which are permuted by the flow and keep their individuality even as $t\to\infty$. 

In \cite{Shn-long} it was proved
that the set $\mathscr{GM}\subset YU(M)$ is nonempty and it is a weak attractor for the Euler flow 
(see also \cite{DoDr}).
(In fact, it is an attractor in a specific sense and it has not been proved that it is an attractor in the usual sense.) 
For some domains (including the periodic strip) 
Bedrossian and Masmoudi \cite{BedMasm2015}
proved that this subset is nontrivial, 
i.e., there exists $u\in YU(M)$ such that $u\notin \mathscr{GM}$. 

Any stationary flow is by definition a GM-flow. 
However, numerical examples show that there exist nonstationary GM flows (at least, on the torus). 
Such a flow comprises several large vortices (blobs) gracefully moving around 
and permanently changing their shapes. These flows  appear to be time-periodic and quasiperiodic; 
it is unclear whether they can be more complex (say, chaotic). 

If $u(x)$ is a stationary flow (i.e. a fixed point of the group $\{S_t\}$) then the level lines of vorticity 
are at the same time the flow lines and hence are real-analytic. 
Futhermore, if $u\in \mathscr{GM}$ is, say, periodic or quasiperiodic, 
our conjecture is that the level lines of vorticity ${\rm curl}\, u(x,t)=\text{const}$ are analytic as well. 
At least, this property is preserved by the Euler evolution. 
Hence, we propose several problems/conjectures.
\begin{problem} 
Prove that for any two-dimensional compact analytic Riemannian manifold $M$ with analytic boundary 
the set $\mathscr{GM}(M)$ of generalized minimal flows is a nonempty and proper subset of 
the Yudovich space $YU(M)$. 
\end{problem} 
\begin{problem} 
Prove that for any GM-flow $u(x,t)$ the level lines of vorticity ${\rm curl}\, u(x,t)=\text{const}$ are real-analytic. 
\end{problem} 

\begin{figure}
\centering
\includegraphics[width=5.5in]{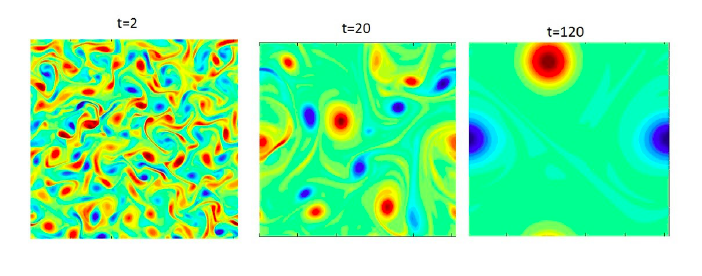}
\vspace{-7pt}
\caption{\small  The appearance of vorticity blobs in the flows on the two-torus.}
\vspace{-10pt}
\label{Fig.2}
\end{figure}
\medskip

One other property of GM-flows is observed in the numerical simulations, cf. Figure~\ref{Fig.2}. The flow domain $M$ 
contains some number of vorticity blobs $B_1, \ldots, B_N$ and a background $B_0$. 
Accurate numerical results hint at the following conjecture.
\begin{problem} 
For any GM-flow ${\rm curl}\, u=\text{const}$ in $B_0$. 
\end{problem} 
Thus, the level sets of ${\rm curl}\, u$ are of two sorts: real-analytic curves and entire domain $B_0$. 

\medskip

\subsection{Are anisotropic function spaces  the future of hydrodynamics?}
The function spaces used in the theory of linear PDEs are less suitable for the nonlinear ones 
and they seem completely inadequate for the Euler equations. One of the most prominent features of the Euler solutions 
is the built-in analyticity of the particle trajectories, flow lines for stationary flows, 
and (conjectured) analyticity of level lines of vorticity of GM-flows. 
\begin{problem} 
Define  anisotropic spaces more appropriate for hydrodynamics with larger smoothness in the direction of flow lines. 
Prove existence and uniqueness theorems for ideal hydrodynamics in those function spaces in any dimension. 
\end{problem} 

The concept of GM-flows appears closely related to the ``ancient flows" for the Navier--Stokes equations 
but, in fact, it is different. Both concepts express the general idea that, in reality, we observe flows that have existed 
since long time ago (one may imagine a river flow). In this regard, it is not natural to consider 
flows starting at some definite time $t_0$ with some (arbitrary) initial velocity $u_0$. 
GM-flows have natural, intrinsic structure (analytic level lines of vorticity) and some natural, finite regularity 
in the transversal direction. 
The question of which functional classes are suitable for their study is an intriguing open problem. 
\begin{problem} 
Determine the true regularity of GM-flows in the framework of anisotropic function spaces. 
\end{problem}


\bigskip

\section{More differential geometry: the $L^2$ exponential map and its singularities} \label{sec:exp}
In classical finite-dimensional Riemannian geometry singular values of the exponential map are the conjugate points. 
The question concerning existence of conjugate points in diffeomorphism groups $\mathscr{D}_\mu^s(M)$ 
and their role in hydrodynamics was posed by Arnold in \cite{arn66}. 
Examples for the two-torus $\mathbb{T}^2$ and the spheres $S^n$ 
were constructed in \cite{mi0, mi} 
with further examples in \cite{sh, pr, pr2, be, wapr, dmsy, tauyon, lty}. 

In infinite dimensions conjugate points are of two types depending on whether 
the derivative of the exponential map fails to be one-to-one (mono-conjugate points) 
or onto (epi-conjugate points). 
Moreover, they can accumulate along finite geodesic segments or have infinite order. 
Such pathological situations can be ruled out in 2D but not in 3D hydrodynamics. 
More precisely, we have the following result from \cite{emp, mp}: 
\begin{thm} 
Let $M$ be an n-dimensional Riemannian manifold (possibly with boundary). 
\begin{enumerate} 
\item[(i)] 
If $n = 2$ then the $L^2$ exponential map is a nonlinear Fredholm map of index zero. 
\item[(ii)] 
If $n \geq 3$ the the Fredholm property of the $L^2$ exponential map fails in general. 
\end{enumerate} 
\end{thm} 
Thus the structure of singularities of the exponential map in 2D hydrodynamics 
resembles that of a smooth map between finite-dimensional manifolds. 
This leads to the following natural question. 
\begin{problem} 
Quantify the failure of the Fredholm property of the $L^2$ exponential map in $\mathscr{D}^s_\mu(M)$ 
for a compact Riemannian 3-manifold $M$. 
\end{problem} 
Examples of three-dimensional manifolds for which the exponential map \eqref{eq:exp} is not Fredholm 
can be found for example in \cite{emp, mp, pr}. 
It is reasonable to expect that 
this failure is borderline 
%
in the following sense. Explicit formulas for the derivative of 
the exponential map of a general right-invariant Sobolev $H^r$ metric derived in \cite{mp} 
decompose it as $A + K_r$ where $A$ is an invertible operator and 
$$
w \mapsto K_r w = (\textbf{1} - \Delta)^{-r/2} P_e ( \iota_w d (\textbf{1}-\Delta)^{r/2} v_0^\flat )^\sharp, 
\qquad 
v_0 \in T_e\mathscr{D}_\mu^s 
$$ 
where $P_e = \textbf{1}-\nabla\Delta^{-1}\mathrm{div}$ 
is the usual Helmholtz--Weyl (or Leray--Hopf) projector onto divergence-free vector fields, 
$\iota_w$ denotes the interior multiplication by a vector field $w$, 
while 
$^\flat$ and $^\sharp$ stand for the standard isomorphisms of the Riemannian metric on $M$. 
When $n \geq 3$ and $r>0$ then the operator $K_r$ turns out to be compact 
and the associated $H^r$ exponential map is Fredholm. 
In the case of the 3D fluids $K_0$ is no longer compact. 
%
In order to measure the deviation of $K_0$ from being a compact operator 
one can, for example, examine its essential spectrum. 
Explicit examples like the rotating solid cylinder or the Taylor--Green vortex, as well as careful numerical experiments, 
may provide some valuable insight. 


Here are several other interesting questions concerning the structure and the role of 
conjugate points in 2D hydrodynamics. 
\begin{problem}\label{prob-cpts} 
\begin{itemize} 
\item[] 
\item[(a)] 
Determine the order of the first conjugate point along any $L^2$ geodesic starting from the identity 
in $\mathscr{D}_\mu^s(M)$. 
\item[(b)] 
Is there a relation between existence of conjugate points in $\mathscr{D}_\mu^s(M)$ 
and Arnold stability criterion for stationary flows in $M$? 
\end{itemize} 
\end{problem} 
%



The answer to $(a)$ may have some bearing on Problem \ref{2ptbvp}. 
For example, if the order of any conjugate point turns out to be always greater than one 
then the fact that the exponential map is Fredholm of index zero would indicate that 
there is only one connected component of the identity over which $\exp_e$ is a covering map, 
see \cite{mp}. 

Regarding $(b)$ recall that according to Arnold's criterion \cite{arn66, MatMet, ak} 
a stationary flow of an ideal fluid is Lyapunov stable if the quadratic form given by 
the second derivative of the kinetic energy restricted to the coadjoint orbits is positive or negative definite. 
It can be shown that for simple domains such as the disk, the annulus and the straight channel 
no steady flows satisfying Arnold's stability possess conjugate points. (We assume here that 
the two-dimensional fluid domain $M$ has a nonempty boundary.) 
It is therefore tempting to expect that this is true more generally. 
See \cite{dmsy} for additional background and \cite{tayo} for recent results in this direction.

\section{Long time behaviour of 2D flows} 

Since existence, uniqueness and regularity of 2D solutions of \eqref{ideal} on the infinite time interval 
is quite well established we can proceed to ask questions concerning the long time behaviour of fluid flows. 

\subsection{Complexity growth for 2D flows}  
Let $M$ be a two-dimensional compact manifold possibly with boundary. 
Recall that the $L^2$ exponential map \eqref{eq:exp} is a local diffeomorphism near the identity $e$ 
in $\mathscr{D}_\mu^s(M)$. 
Fix $\epsilon >0$ and let 
$\mathcal{U}_\epsilon = \exp_e(B_\epsilon)$ 
where 
$B_\epsilon = \{ v \in T_e\mathscr{D}_\mu^s~|~ \| v \|_{H^s} < \epsilon \}$ 
is an open $H^s$ ball of radius $\epsilon$. 
Any diffeomorphism in $\mathscr{D}_\mu^s(M)$ can be represented as a product 
$\eta = \eta_1 \circ \cdots \circ \eta_N$ 
of a finite number of elements from $\mathcal{U}_\epsilon$, 
see \cite{luk0, luk2}. 
Let $\mathcal{C}_\epsilon(\eta)$ denote the minimal number of factors in this representation 
and let 
$\mathcal{C}(\eta) = \limsup_{\epsilon \to 0} \big( \epsilon \, \mathcal{C}_\epsilon(\eta) \big)$ 
be the \emph{absolute complexity} of $\eta$. 
\begin{problem} 
Show that for any geodesic $\gamma(t)$ in $\mathscr{D}_\mu^s(M)$ its absolute complexity 
is exponentially bounded above:
$$
\mathcal{C}(\gamma(t)) \lesssim e^{t\|\dot{\gamma}_0\|_{H^s}}\,.
$$ 
\end{problem} 
This estimate would imply the well known double exponential estimate 
for solutions of the Euler equations \eqref{ideal}, namely 
$\|v(t)\|_{H^s} \lesssim e^{C_1e^{tC_2}}$. 
However complexity of a flow as defined above has a broader sense than its regularity. 
\begin{problem} 
Show that for a typical geodesic $\gamma(t)$ in $\mathscr{D}_\mu^s(M)$ 
we have 
$\mathcal{C}(\gamma(t)) \simeq t$. 
\end{problem} 
Roughly, we say that a family of geodesics starting from the identity in $\mathscr{D}_\mu^s(M)$ 
is \emph{typical} if the complement of the corresponding set of initial velocities in $T_e\mathscr{D}_\mu^s$ 
has infinite codimension. 
Examples of ``typical" $L^2$ geodesics whose complexity grows linearly in time 
are the one-parameter subgroups of $\mathscr{D}_\mu^s(M)$ 
and (possibly) the quasi-periodic solutions of \eqref{ideal}. 
%

\subsection{Aging of the fluid, irreversibility and Lyapunov functions.} 
Consider an arbitrary solution $u=u(t,x)$ of the Euler equations in $M$. 
Given any two time instants can one decide which velocity field $u(t_1,x)$ or $u(t_2,x)$ 
corresponds to an earlier time instant? 
In other words, is it possible to determine the aging of the fluid from its velocity field? 
Numerical experiments suggest that fluids ``age" with time: 
starting with a smooth initial velocity $u_0(x)$ the corresponding solution becomes ``wrinkled" 
in that its derivatives generally grow.\footnote{V.~Yudovich referred to this phenomenon as 
``regularity deterioration".} 
\begin{problem}[Aging Problem] 
Is it possible to quantify the ``aging" property of the fluid? 
\end{problem} 
Perhaps the best way is to find a Lyapunov function $L$ defined and continuous on the space of fluid velocities 
in $T_e\mathscr{D}_\mu^s$ such that 
$\tfrac{d}{dt}L(u(t,\cdot)) \geq 0$ 
where $u$ is a solution of \eqref{ideal} 
with equality holding on a ``slim" subset of $T_e\mathscr{D}_\mu^s$ (say, of infinite codimension). 

The first Lyapunov function in this context was constructed by Yudovich in 1970s 
for flows with a rectilinear streamline 
(e.g., flows in domains whose boundary contains a straight line segment), 
see \cite{yu1, yu2}. 
The construction was subsequently generalized for arbitrary bounded domains in \cite{msy}. 
It implies ``regularity deterioration" at least on the boundary. 

It is natural to expect that there are also Lyapunov functions which are supported inside the fluid domain. 
Examples describing evolution of weak singularities in the Lagrangian flow 
were found in \cite{sh1}. 
The fact that these singularities become gradually ``sharper" as the fluid evolves 
suggests the same deterioration phenomenon stressed by Yudovich. 

Consider an ideal fluid in a periodic channel $M=\mathbb T \times (0,a)$. 
Assume that its initial velocity $u_0$ is $C^1$ close to that of a plane-parallel flow 
whose velocity profile $v=v(x_2)$ satisfies 
$v'>0$ and $v''>0$ in $(0,a)$. 
Suppose that the level lines of $\omega_0 = \mathrm{curl}\,u_0$ satisfy $\nabla\omega_0 \neq 0$ in $M$ 
and that one of the lines $\omega_0 = \text{const}$ has a ``kink". 
\begin{problem} 
Show that the ``kink" does not disappear as the fluid evolves in time. 
\end{problem} 
%

\section{Entropy and fluids} 
%
\subsection{Entropy of a set and entropy of a measure} 
In general terms, entropy is a measure of diversity of some ensemble. 
For a finite set $S$ with $N$ elements of equal ``weights" (i.e., equivalent in some respect) 
the entropy $H(S)$ is equal to $\log_2{N}$. 
If the elements $s_i \in S$ $(i=1, \dots N)$ have different weights $w_i$ 
then we define the entropy of the weighted finite set $S$ to be 
$H(S) = - \sum_{i=1}^N w_i \log_2{w_i}$. 
For example, if the whole mass (assumed to be 1) is concentrated at some $s_k$ then $H(S)=0$; 
otherwise $H(S)>0$ 
with maximum value $\log_2{N}$ if $w_i = 1/N$ for all $i = 1, \dots N$. 

If $S$ is an infinite set then the definition of $H(S)$ is not so clear. 
It is based on approximations of the set $S$ by finite sets 
and of the weight (i.e., the probability measure) $\mu$ on $S$ by some discrete weights. 

Suppose that $S$ is a compact subset of a complete metric space $X$. 
In the absence of a measure on $S$ we can define the Kolmogorov $\epsilon$-entropy of the set $S$ as  
$H_\epsilon(S) = \log_2{N_\epsilon}$ 
where $N_\epsilon$ is the cardinality of the minimal $\epsilon$-net, 
i.e., the minimal number of $\epsilon$-balls in $X$ covering $S$. 
(In fact, $H_\epsilon(S)$ is usually defined as an equivalence class of such functions as $\epsilon \to 0$.) 
If $S$ is equipped with a measure then we define an analogue of a weighted finite set as above 
and we can try to define a suitable analogue of the $\epsilon$-entropy.\footnote{It is not clear if such an object 
has anything in common with the entropy of an invariant measure as typically defined in the theory of dynamical systems.} 
One option is to define it as the $\epsilon$-entropy of the support of the given probability measure $\mu$. 
This, however, would give only an upper bound, just like for weighted finite sets. 
Another option is to define 
$$ 
H_{\epsilon,\delta}(\mu) 
= 
\inf{\big\{ 
H_\epsilon(Y)~|~ \text{$Y \subset X$ is compact and} \; \mu(X\setminus Y) \leq \delta 
\big\}}. 
$$ 
What is the relation of $H_{\epsilon,\delta}(\mu)$ to the entropy of a weighted finite set?

The definition of entropy given above depends on a pre-existing measure on $X$. 
If $X$ is a finite set then we can take $\mu$ to be a counting measure; 
if $X$ is a phase space of classical mechanics then $\mu$ could be a Liouville measure. 
But on a general metric space no such choices are available. 
On the other hand, if $X$ is also a vector space then we can cover its compact subsets 
by congruent sets other than metric balls. For example, we can use cylindrical domains 
with finite-dimensional base. 
Suppose that $X$ is a Hilbert space with coordinates $x_1, x_2, \dots$ and let $X_n$ be the subspace 
defined by $x_i=0, i>n$. 
Subdivide $X_n$ into cubes $C_j$ of side length $\epsilon>0$ and let $K_j = \pi_n^{-1}(C_j)$ 
where $\pi_n$ is the orthogonal projection onto $X_n$. 
Given a compact set $S \subset X$ for any $n$ let $\tilde{N}(n,\epsilon)$ be the number of cubes 
$C_j \subset X_n$ having nonempty intersection with $\pi_n(S)$. 
Let $\tilde{N}(\epsilon)=\sup_n{\tilde{N}(n,\epsilon)}<\infty$ and now define 
$\tilde{H}_\epsilon(S) = \log_2{\tilde{N}(\epsilon)}$. 
Note the similarity of this function to $H_\epsilon(S)$ defined previously. 

Finally, consider an analogue of the $\epsilon,\delta$-entropy for a probability measure $\mu$. 
Let $X_n$ be a finite-dimensional subspace of $X$ as before subdivided into $\epsilon$-cubes $C_j$ 
and let $K_j = \pi_n^{-1}(C_j)$. 
Set $H_{\epsilon,n}(\mu) = \sum_j \mu(K_j) \log_2{\mu(K_j)}$ 
and observe that for any fixed $\epsilon$ this quantity is bounded uniformly in $n$. 
Define the $\epsilon$-entropy of $\mu$ to be 
$H_\epsilon(\mu) = \sup_n{ H_{\epsilon,n}(\mu)}$. 
\begin{problem} 
Investigate  properties of the entropy functions $H_\epsilon, H_{\epsilon,\delta}, \tilde{H}_\epsilon$ 
and $\tilde{H}_{\epsilon,\delta}$ in this section and explain relations between them. 
\end{problem} 

As an example, if $M = Q^d$ is the unit cube in $\mathbb{R}^d$ 
and $\mu=\mu^d$ is the $d$-dimensional Lebesgue measure then 
$H_\epsilon(\mu) = \epsilon^{-d} \epsilon^d \log_2{\epsilon^d} = d \log_2{\epsilon}$. 
The result is the same if we view the cube as a subset $Q^d \subset X_d$ supporting $\mu^d$ 
and compute its $\epsilon$-entropy in the whole space $X$. 
But what will happen if we consider $Q^d$ as a subset of $X_n$ for some $n>d$ 
and rotate it so that it is no longer a coordinate cube? 
In this case the sum becomes 
$\sum_j \mu^d(C_j) \log_2{\mu^d(C_j)} \sim \log_2{\epsilon^d} + { o}(\log_2{\epsilon})$ 
so that the principal asymptotic does not change under rotations of $M$. 
The same can be said about other deformations: 
the $\epsilon$-entropy behaves like $d\log_2{\epsilon} \cdot \mu^d(X)$.

\subsection{Entropy decrease for the Euler flow} 
Computer experiments show that the velocity field of an ideal 2D fluid behaves similarly 
for all initial conditions. 
The outcome is a small collection of moving vortices forming a hierarchical structure of 
``islands", ``lakes", ``satellites", ``archipelagoes", etc. 
that are not mixing but instead preserve their individuality. 
This scenario looks quite strange from the physical viewpoint. 
It is obvious that the diversity of the initial conditions is much higher than that of the outcomes. 
The quantitative measure of the diversity in this case is the $\epsilon$-entropy. 

The natural (physical) phase space here is 
$V^0 = \{ u \in L^2(M,\mathbb{R}^2)~|~ \mathrm{div}\, u = 0, \; u\,\|\,{\partial M} \}$. 
Consider an initial velocity ensemble, i.e., a probability measure $\mu_0$ in $V^0$. 
In fact, the initial velocity should be more regular than merely $L^2$ 
so that the Cauchy problem is correctly posed. 
A good example is the Yudovich space 
$Y = \{ u \in V^0~|~ \mathrm{curl}\, u \in L^\infty \}$, 
but we can also start with an initial velocity in $V^s=T_e\mathscr{D}_\mu^s$ for some $s>2$. 
The exact value of the Sobolev index $s$ is not essential because in the long run 
the flow will approach some asymptotic regime possessing a ``natural" though presently unknown regularity. 
This asymptotic flow will belong to the Yudovich space and may be more regular 
but it is not clear if this regularity can be captured by some nice function space 
(e.g., it may be an element of a Frechet space and there may be no reasonable choice of 
a Banach space for this purpose). 

Suppose that $\mu_0$ is a compactly supported measure in $V^0$. 
We can define its $\epsilon$-entropy $H_\epsilon(\mu_0)$ 
(or at least its principal asymptotic as $\epsilon \to 0$). 
Let $\mu_t$ be the measure at time $t$ transported by the Euler flow. 
The phenomenon discussed here can be described by the inequality 
$\liminf_{t \to \infty} H_\epsilon(\mu_t) \preceq H_\epsilon(\mu_0)$ 
in the sense that for a limiting measure $\mu_\infty$ we have 
$H_\epsilon(\mu_\infty) = \mathcal{O}(H_\epsilon(\mu_0))$ as $t \to \infty$ 
and for some measures $\mu_0$ the ``big Oh" is replaced by the ``little oh". 

Such a result would make a physicist uneasy because it looks like a violation of the Liouville theorem. 
However, this is not a real violation since the true phase space of the fluid includes not only velocities 
of  all the fluid particles but also their positions. 
Therefore, the elements of the space $V$ represent merely ``half" of the phase space coordinates 
with the other ``half" represented by the flow map in $\mathscr{D}_\mu$. 
\begin{problem} 
Explain the phenomenon of entropy decrease for solutions of the Euler equations \eqref{ideal}. 
\end{problem} 

\bigskip

\section{Hamiltonian properties of the Euler equation}

The Riemannian geometric approach to hydrodynamics has a Hamiltonian reformulation, 
see e.g. \cite{ArnHamiltonian, ak}. 
Namely, consider again the group of smooth volume-preserving diffeomorphisms  $\mathscr{D}_\mu(M)$
and denote its Lie algebra of divergence-free vector fields by $\mathfrak g = \Vect_\mu(M)$. 
The regular dual space $\mathfrak g^*$ of $\mathfrak{g}$ can be naturally identified with 
the space   of 1-forms modulo differentials of functions  on $M$, i.e.  
with the space of cosets  $\mathfrak g^*=\Omega^1(M)/d\Omega^0(M)$. 
The inertia operator $A: \mathfrak{g} \to \mathfrak{g}^\ast$ relies on the choice of metric on $M$ 
and to a divergence-free vector field $v$ it associates the coset $[v^\flat]$ of the 1-form $u=v^\flat$ 
related by means of the metric.
Then the hydrodynamic Euler equation \eqref{ideal} can be written as an evolution of 1-forms
$$
\partial_t u+L_v u=-df
$$
for a certain time-dependent function $f$ on $M$, or as an evolution of cosets of 1-forms:
$$
\partial_t[u]=-L_v[u],
$$
where $u=v^\flat \in\Omega^1(M)$ and $[u]\in\Omega^1(M)/d\Omega^0(M)$.
This is a Hamiltonian equation  with respect to the Lie-Poisson structure on $\mathfrak g^*$ and 
with the Hamiltonian function given by the fluid's kinetic energy 
$E(v)=\frac{1}{2} \langle Av, v\rangle = \frac{1}{2} \| v\|_{L^2(M)}^2$. 
This way,  the equations of an ideal fluid dynamics in any dimension form a Hamiltonian system. 

\subsection{Nonintegrability}
In the two-dimensional case, besides the kinetic energy, this system has infinitely many enstrophy invariants. 
These invariants are Casimir functionals: they do not depend on the metric on $M$, but specify a coadjoint orbit 
(i.e. a particular set of isovortical vector fields) on which the Euler evolution takes place but say nothing 
about the dynamics on the orbit itself.
\begin{problem}
Prove the non-integrability of the 2D Euler equation.
\end{problem}

It is worth to emphasize that, while there is a number of papers related to infinitely many conserved quantities 
or the Lax form of the Euler equation in 2D (see, e.g., \cite{fr-vishik, li-yur}),
these features cannot be regarded as good indicators of integrability.
The existence of a Lax pair is a property of all Euler--Poincar\'e equations, since the latter  
are Hamiltonian on the dual space to a Lie algebra with respect to the Lie--Poisson structure 
and hence are given by the coadjoint operator, ``mimicking the commutator" in the Lie algebra. 
Typically, in order to prove algebraic-geometric integrability of a Lax equation
one needs to present a Lax form nontrivially depending on a spectral parameter. 
However, no such form has been found for the 2D Euler equations. 

On the other hand, one could try 
to prove non-integrability of this infinite-dimensional system with the help of the methods 
used to  show non-integrability of finite-dimensional Hamiltonian systems.
This could invoke, for instance, the methods related to nontrivial monodromies for periodic orbits 
 \cite{ziglin, mor-ramis}, or Melnikov integrals for bifurcations of saddle separatrices, or
 one might be able to proceed by means of a ``local" analysis in the vicinity 
 of a steady solution. 
Since there are several (not necessarily equivalent) definitions of integrability in infinite dimensions, 
the above problem would be to show that the 2D Euler equations fail to satisfy at least one of these 
integrability definitions. 

Note that the set of all Casimirs  has been recently fully described  in \cite{IK17, IKM16} 
for two-dimensional surfaces $M$ without boundary for the groups of symplectic and Hamiltonian diffeomorphisms. 
The two-dimensional boundary case was settled in \cite{Kiril}.
\medskip

\subsection{Finite-dimensional approximations}
It is also worth pointing out that 
there are various approximations of the 2D Euler equations on the plane, the 2-torus or the 2-sphere 
by finite-dimensional  Hamiltonian systems. Two of the best known approximations are 
$(i)$ by ${\rm SU}(N)$-algebras, whose structure constants converge to those of the algebra 
$\rm{Vect}_\mu(\mathbb{T}^2)$ \cite{zei}
and 
$(ii)$  by a system of $N$ point vortices on $\R^2$, as $N\to\infty$ \cite{mar-pul}.
In a recent paper \cite{ModViv21} it was shown how both of these approximations 
can be unified within a quantization approach to the 2D hydrodynamics. 

\medskip

Dynamics of $N$ point vortices on the plane has been studied since the time of Helmholtz and Kirchhoff \cite{mar-pul} 
and is of particular interest; see e.g. \cite{Newton01}. 
For instance, on the plane or the sphere the corresponding finite-dimensional Hamiltonian systems 
turn out to be integrable for $N\le 3$, while for the torus the system is integrable only for $N\le 2$.
This is related to the fact that the corresponding Kirchhoff equations for point vortices on $\mathbb R^2$ and $S^2$ are invariant with respect to the three-dimensional isometry groups $E(3)$ and $SO(3)$, while for $\mathbb{T}^2$ this group of isometries is isomorphic to $\mathbb R^2$ and hence two-dimensional. 
While integrability  still holds for 4 vortices at zero total vorticity on the sphere $S^2$
and at zero total vorticity and momentum on the plane $\mathbb R^2$, for $N\ge 4$ point vortices of generic strengths the motion becomes non-integrable, see \cite{KozlovTre, ziglin82}.
\medskip

In \cite{ModViv19, ModViv21, Shn-long} the authors observed the following intriguing phenomenon:  
on the two-dimensional sphere and torus the 2D Euler motion for a smooth and sufficiently general  initial vorticity 
after some time (and with a very small numerical viscosity) leads to merging of smaller vortex formations of the same sign 
into larger blobs, cf. Figure~\ref{Fig.2}. Surprisingly, in those numerical
simulations  it was recovered that this clustering continues until the blob dynamics,
approximated by point vortex motion, ``becomes integrable".  Namely, such clustering leads to an integrable dynamics of 2 vortex blobs on the torus, 
 3  vortex blobs for general angular momentum on the sphere, 
 and 4 vortex blobs  on the sphere provided that  the initial angular momentum was zero, thus 
 exactly recovering integrable cases for $N=2,3,$ and $4$ point vortices on $\mathbb{T}^2, \mathbb R^2$ and $S^2$ discussed above,
 see \cite{ModViv19, ModViv21}.

\begin{problem} 
Justify the phenomenon observed in {\rm \cite{ModViv19, ModViv21}} that integrable cases of point vortices 
seem to be attractors for 2D Euler flows with generic smooth initial vorticity in 2D (for small numerical viscosity). 
Design  a model of dissipation in the 2D Euler equation which would produce such an integrable dynamics at large times.
\end{problem}

 \medskip

Historically, studies of point vortices include constructions of explicit solutions, 
conditions for the presence and absence of a collapse, 
description of relative equilibria, bifurcation of solutions for $N=2 $ and $N=3$ point vortices, 
see \cite{Aref, Newton01}. 
A nice collection of integrable motions on the sphere and torus can be found in \cite{ModViv20, Newton01}.
For manifolds with boundary there is a broader variety of motions. 
For instance, the cusp motion of a pair of point vortices on a half-plane is related to the golden ratio \cite{KWang}. 
The motion of 3 point vortices on a half-plane is already non-integrable \cite{Cheng}, 
as well as, apparently, the motion of two point vortices in the quadrant. 
The motion of point vortices on the half-plane and the
quadrant for the lake equation is related to the 
motion of vortex rings, membranes and to the more general binormal equation 
\cite{Jer02, JerSm15, Cheng}.
\begin{problem} 
Study in detail the motion and bifurcations of a small number of vortices on various manifolds: 
half- and quarter-plane, hemi- and quarter-sphere, disk, torus, cylinder and half-cylinder, etc. 
\end{problem}

Of particular interest is the study of point vortices on {\it non-orientable manifolds}, 
which started only recently, see e.g. \cite{balab2, balab1, Mobius, enric} for a description  of  the motion and bifurcations of a small number of vortices 
on such non-orientable surfaces as the M\"obius band, projective plane and the Klein bottle. 
\begin{problem} 
Study the first non-integrable cases for a small number of point vortices on non-orientable surfaces.  Describe the full set of Casimirs and finite-dimensional approximations for the group of area-preserving diffeomorphisms of non-orientable surfaces, analogous to the ones in the orientable cases. 
 \end{problem}

\bigskip

\section{Dynamical properties of the Euler equation: wandering solutions, chaos, non-mixing and KAM.} 
Consider a finite-dimensional model of Euler hydrodynamics: the Euler--Poincar\'e equation 
on a finite-dimensional Lie group corresponding to a positive-definite energy form. 
This system is Hamiltonian on a coadjoint orbit and satisfies the condition of the Poincar\'e recurrence theorem. 
Indeed,  by fixing the energy level one confines the dynamics to the compact set (even for a non-compact group) 
which is the intersection of the energy level and the orbit. 
The dynamics preserves the volume form on this intersection (see \cite{ak}) 
and hence yields to Poincar\'e's recurrence. Therefore every point of the orbit in the course of evolution 
returns arbitrarily close to its initial position after arbitrarily large time. 

This is not the case for a general infinite-dimensional dynamical system 
and, in particular, for the Eulerian hydrodynamics. 
Nadirashvili showed that the Euler equation of a 2D fluid has wandering solutions: 
there is an initial condition of a fluid in a 2D annulus whose neighbourhood never returns sufficiently closely 
to the initial condition after a certain time \cite{Nadir}.
A 3D analogue of that result is unknown. 
\begin{problem} 
Prove that the 3D Euler equation has  wandering solutions. 
\end{problem} 

\medskip

The only type of results in this direction are the non-transitivity and non-mixing properties of 
the 3D Euler equations proved using the ideas of the KAM theory; see  \cite{KKP14}.
It turns out that the dynamical system defined by the hydrodynamical Euler equation 
on any closed Riemannian 3-manifold $M$ is not mixing in the $C^k$ topology (for $k > 4$ and non-integer) 
for any prescribed value of helicity and sufficiently large energy. 
Furthermore, this non-mixing property of the flow of the 3D Euler equation has a local nature: 
in any neighbourhood of a ``typical'' steady solution on $\mathbb S^3$ there is a generic set of initial conditions 
such that the corresponding Euler flows will never enter a vicinity (in the $C^k$ norm for any non-integer $k>10$) 
of the original steady flow; see \cite{KKP19}. 
 
Along the way one constructs a family of functionals on the space of divergence-free $C^1$ vector fields on $M$ 
which are integrals of motion of the 3D Euler equation: 
given a vector field these functionals measure the part of the manifold $M$ foliated by ergodic invariant tori 
of fixed isotopy types. The KAM theory allows one to establish certain continuity properties of these functionals 
in the $C^k$-topology and to get a lower bound on the $C^k$-distance between a divergence-free field 
(in particular, a steady solution) and a trajectory of the Euler flow. This way one obtains an obstruction 
for the mixing under the Euler flow of $C^k$-neighbourhoods of divergence-free vector fields on $M$. 
The local version of non-mixing is based on a similar KAM-type argument to generate knotted invariant tori 
from elliptic orbits in nondegenerate steady Euler  flows. 
\begin{problem} 
Relax the restrictions on the smoothness index $k$ of the $C^k$ spaces 
(related to application of the KAM) and prove non-transitivity and non-mixing properties 
of the 3D Euler in full generality. 
\end{problem}

\smallskip

It turns out that the Euler equations on higher-dimensional Riemannian manifolds 
possess a kind of universal embedding property, somewhat similar to the theorems of Whitney and Nash 
on embeddings of manifolds as submanifolds in higher-dimensional Euclidean spaces. 

Namely, it was shown in \cite{Tao} that a certain large class of 
finite-dimensional quadratic dynamical systems in $\R^d$ can be realized as subsystems 
of the hydrodynamical Euler equation on the manifold ${\rm SO}(d)\times \mathbb T^{d+1}$ 
with a certain metric depending on the original system. 
Subsequently, Torres de Lizaur in \cite{Torres} proved that such a realization is possible 
for any dynamical system of a finite-dimensional manifold or for its approximation. 
This way essentially any finite-dimensional dynamical system or its approximation to an arbitrary degree 
can be embedded as an invariant (tiny) subsystem in a higher-dimensional Euler equation for a certain metric. 

The construction in \cite{Torres}, which has already found other applications, goes as follows: 
for a given finite-dimensional dynamical system one first embeds it via Whitney to a system 
on a submanifold inside a higher-dimensional torus, extends it hyperbolically to a
dynamical system in the torus, and then writes it via the smooth
vector field (represented by a Fourier series) in that torus. Then one
truncates it to a Fourier polynomial (this is where the approximation
with an arbitrary precision takes place).  Finally, one observes that
a Fourier polynomial vector field $v(x)=\sum_l c_l \,e^{\mathrm{i}lx} \,{\partial}/{\partial x} $ rewritten in trigonometric
coordinates $p_l:=e^{\mathrm{i}lx}$ becomes quadratic:  
$v(p)=\mathrm{i}\sum_{k,l}  k \,c_l\, p_k p_l \,{\partial}/{\partial p_k}$.
After that one employs Tao's embedding \cite{Tao} of quadratic systems to the higher-dimensional Euler equations. 
The examples include such structurally stable systems exhibiting chaos as the  ABC flows inside the higher Euler phase space.

One should note that the dynamics of the Euler equation  outside of this tiny submanifold is not controlled and, in principle, could be rather regular. For instance, 
one might have dynamical systems with a very regular 
behaviour almost everywhere, but with some chaotic behaviour on a very tiny submanifold.

\begin{problem}
Consider an integrable system on a compact $2n$-dimensional manifold, which 
has $n$ first integrals in involution, functionally independent almost everywhere.
How wild could  such a system be on a singular submanifold, where the integrals become dependent? What are constraints dictated by the integrability? Could one observe a chaos on some tiny submanifold of an integrable system?
\end{problem}

\section{Steady Euler flows} 
Steady Euler flows in a domain $M\subset \R^d$ are defined by the equation $v\cdot\nabla v=-\nabla p$
along with the divergence-free restriction $\mathrm{div}\, v=0$
and the condition   $v||\partial M$ of  tangency to the boundary. 
While it is easy to construct a steady 2D flow with compact support in $\R^2$ 
(take a radial stream function with compact support) it is a notoriously difficult task to perform this feat in 3D.

 An explicit recent example of a smooth steady incompressible Euler flow in $\R^3$ with compact support 
 was given in \cite{ga}, see also a more general approach of \cite{co19}. 
 In this type of solutions the pressure and Bernoulli function are dependent. 
 
Arnold pointed out the remarkable topology of steady 3D fields: 
for an analytic steady field, not everywhere collinear with its curl, the flow domain is almost everywhere 
fibered into invariant tori and annuli, see \cite{arn66, ak}. If the steady field is everywhere collinear with its vorticity 
but the proportionality coefficient is a generic function, then the domain is still fibered in a similar way. 
If the steady field is Beltrami, i.e., it is an eigenfield for the curl operator, 
$\mathrm{curl}\,v=\lambda\,v$, then its topology can be very intricate. 
Hence a paradox arises: 
a generic steady field has a very regular topology, while a sufficiently chaotic field must necessarily 
be an eigenfield for the curl operator. 
\begin{problem} 
Explain this paradox: what is a typical steady field and what genericity notion is natural for steady fields? 
\end{problem} 

It is worth mentioning that the notion of a typical object in fluid dynamics might be quite different from the standard one. 
For instance, as discussed in Section \ref{sect:analytic} stream functions of steady 2D flows always have analytical levels 
even if they have only finite smoothness across the levels.

In \cite{myz} the authors considered slightly compressible 3D vector fields and their incompressible limit 
to explain the above paradox. 
The paper \cite{EP14} also sheds more light on this problem: the levels of the Bernoulli function cannot be spheres. 

\bigskip

It turned out that the topology of Beltrami fields can be arbitrarily complicated. 
In \cite{EP11} it was established that for any finite link $L\subset \R^3$ and any nonzero real number $\lambda$ 
one can deform the link $L$ by a $C^\infty$ diffeomorphism of $\R^3$, arbitrarily close to the identity in any $C^m$ norm, 
such that the image of the link becomes a set of vortex lines of a Beltrami field $v$ with the eigenvalue $\lambda$ in $\R^3$, 
${\rm curl }\, v =\lambda v$  in $\R^3$ and, moreover, $v$ falls off at infinity as $|x|^{-1}$. 
In \cite{tubes} a similar result was proved for the existence of a finite collection of 
toroidal knotted or linked stream/vortex tubes in $\R^3$. 
The boundaries of such tubes are structurally stable invariant tori for a Beltrami field with a quasiperiodic flow on them. 
Furthermore, there are Beltrami fields with invariant tori of arbitrary topology that enclose regions 
with any prescribed number of hyperbolic periodic orbits, see \cite{arxiv, EPT}. 
In this series of papers  Encisco and Peralta-Salas with coauthors established 
other interesting topological properties of Beltrami fields on various manifolds. 
We refer to these papers for various  open problems related to this topic. 
 
 \bigskip
 
Another active direction of research is related to the interaction of 
topological and metric properties of divergence-free vector fields.
By \emph{topological} we mean those properties that are defined using the volume form only, 
e.g. average linking of the field trajectories, which is given by the field's helicity. For an exact divergence-free 
 vector field $ u $ on a three-dimensional manifold $M$ with a volume form $\mu$ its helicity (or asymptotic Hopf) invariant is 
 $$H( u )=\int_M \omega\wedge d^{-1}\omega =\int_M (  u , \text{curl}^{-1} u )  \, \mu\,,
 $$
 where $\omega := \iota_u \mu$ (interior product) is the 2-form on $M$ whose kernel field is $u$, see \cite{Mo69, Mo85, Mo21}. 
 (Here the second expression for helicity, convenient in explicit computations, 
 relies on a choice of a Riemannian metric on $M$, while the first one shows that
 helicity does not depend on that choice.) Actually, the helicity was shown to be the only topological invariant in a large class of functionals under the action of the group of volume-preserving diffeomorphisms  \cite{PNAS}.
(One should also mention that for velocity fields $v$ that are solutions of the Euler equation, their helicity 
is defined as the helicity invariant of the corresponding vorticity field $ u :=\text{curl}\,v$, hence in terms of the 
velocity the corresponding expression is $H( \text{curl}\,v )=\int_M ( \text{curl}\,v, v)  \, \mu$.)
\smallskip

While topological properties of the fields require only fixing a volume form, 
their {\it geometric} properties require a Riemannian metric to define them. 
An example of the latter is the $L^2$-energy of the field 
$E( u ):= \frac{1}{2} \| u \|_{L^2(M)}^2 = \frac 12 \int_M (  u ,  u )  \, \mu$. 

\medskip

The inequality ``helicity bounds energy",
$$ 
E( u )\ge {\rm const}\cdot H( u ) \,,
$$ 
means that  nontrivial average linking of the (magnetic) field's trajectories prevents its energy from complete dissipation 
via volume-preserving diffeomorphisms, see \cite{Arnold73, Mo85}. (This process is often called magnetic relaxation.)
This inequality can be proven by noticing that the operator $\text{curl}^{-1}$ on a compact manifold or domain $M$ 
has bounded spectrum and by applying the Poincar\'e inequality; see  \cite{Arnold73, ak}. 
Geometrically one can visualize this inequality for a vector field confined to a pair of simply linked solid tori. 
To minimize the energy of this linkage one needs to shorten trajectories of the field. 
On the other hand, due to the incompressibility property the shrinking of trajectories in one of the tori 
leads to stretching of the trajectories in the other. 
It is a particularly challenging problem to describe and analyze the process of  magnetic relaxation 
to an equilibrium, possibly nonsmooth; see \cite{Mo85, Mo21, Bee}. 

Note that the ``helicity--energy" inequality is far from being sharp: 
helicity $H(u)$ could be zero, while the field $u$ could possess nontrivially linked tori with opposite linkings 
or a higher order nontrivial linking. 
For knots there is a hierarchy by the Milnor and Massey numbers: 
once the preceding invariants are equal to zero, the invariants of the next level are well defined 
and distinguish the corresponding knots and links. 
\begin{problem}[\cite{ak}] 
Find a sequence of higher helicity invariants for vector fields so that, given  a field, if all the previous invariants are equal to zero for it, then the first nonzero invariant bounds the field's energy from below. 
\end{problem} 

We refer to the book \cite{ak} and its second edition for the discussion of open problems 
and a large bibliography on the subject, see also \cite{fr-he, TuGl13, LauStr, Mo21}. 

\bigskip

\section{Singular vorticities in the Euler equation} \label{sect:vortex}
The localized induction approximation (LIA) procedure applied to the 3D Euler equation in vorticity form 
gives the vortex filament (or binormal) equation: 
$$ 
\partial_t\gamma=\gamma'\times\gamma'' 
$$ 
for the vorticity $\delta_\gamma$ supported on a curve $\gamma\subset \R^3$. 
Similarly, for the vorticity 2-form $\delta_P$ supported on a vortex membrane, 
a submanifold $P^{n-2}\subset \R^n$ of codimension 2, the LIA equation turns out to be 
$$ 
\partial_t q=J({\bf MC}_P(q)), 
$$ 
where ${\bf MC}_P(q)$ is the vector of the mean curvature of the membrane $P$ at the point $q\in P$ 
and $J$ rotates by $\pi/2$ this vector in the normal plane to $N_qP$ to $P$, 
see \cite{Jer02, JerSm15a, HaVi, Shashi, BK-membranes}. 

The binormal equation is known to be equivalent to a 1D compressible fluid equation 
and to the nonlinear Schr\"odinger equation (NLS) in 1D via the Hasimoto transform. 
It also gives singular solutions of the Gross--Pitaevsky (or the 3D NLS) equation \cite{Jer02, JerSm15}. 
\begin{problem} 
Find a direct link from 3D NLS equation to 1D NLS equation as a reduction to singular solutions, 
rather than going through the LIA procedure. 
A similar question arises for the compressible Euler equations as a reduction 
from 3D to singular solutions supported on 1D submanifolds. 
\end{problem} 

\medskip

Another interesting case of singular solutions is that of the vorticity supported on a hypersurface, called a vortex sheet. 
One can introduce a symplectic structure on vortex sheets similar to the Marsden-Weinstein symplectic structure 
on membranes; see \cite{BK-membranes} for the corresponding Hamiltonian formalism.

However, it is more natural to describe the motion of vortex sheets by means of 
a variational principle {\it \`a la} Arnold, 
albeit for a different object, as geodesics on an infinite-dimensional  Lie groupoid,  
see Section \ref{s:groupoids}. 
Using the corresponding vortex sheet groupoid instead of the  group of volume-preserving diffeomorphisms in Arnold's framework, one obtains a geometric interpretation for discontinuous fluid flows, as well as their Hamiltonian description 
on the corresponding dual Lie algebroid, see \cite{IK18}. 
It turns out that vortex sheet type solutions of the Euler equation are precisely the geodesics of 
an $L^2$-type right-invariant (source-wise) metric on the Lie groupoid of discontinuous volume-preserving diffeomorphisms. 
The geodesics on the groupoid turn out to be weak solutions of the Euler equation with vortex sheet initial data \cite{IK18}. 

\medskip

Geometric description of vortex sheets leads to an interesting non-local metric of hydrodynamical pedigree 
on shape spaces: 
it is an $H^{-1/2}$-metric on closed hypersurfaces bounding the same volume 
(or equivalently, on constant densities inside those hypersurfaces), 
see \cite{IK18}.  Such a metric is constructed with the help of the  Neumann-to-Dirichlet operators 
and it is always nondegenerate since it is bounded below by the Kantorovich--Wasserstein distance. 
Geodesics with respect to this metric describe motions of potential fluid flows with vortex sheets 
and fluids with free dynamic boundary, cf. \cite{IK18, Marsden, Loe}. 
\begin{problem} 
Describe the differential geometry of shape spaces equipped with such $H^{-1/2}$-metrics obtained as metrics on dynamic boundaries or vortex sheets.
\end{problem}

\bigskip

\section{The compressible Euler equation and the NLS equation}
There is a well-known relation between the NLS equation and (quantum) compressible fluids in any dimension. 
In 1927 E.~Madelung \cite{Mad} gave a hydrodynamic formulation of the Schr\"odinger equation. 
For a pair of real-valued functions $\rho$ and $\theta$ on an $n$-dimensional manifold $M$ 
(with $\rho>0$) the Madelung transform is the mapping $\Phi :(\rho,\theta)\mapsto \psi$
given by $\psi = \sqrt{\rho e^{\mathrm{i}\theta}}: M \to \mathbb C$. 

The Madelung transform maps the system of equations for a barotropic-type fluid 
to the Schr\"{o}dinger equation. 
Namely, let the function (or density) $\rho$ and the potential velocity field $v=\nabla \theta$ 
satisfy the following barotropic-type fluid equations:
\begin{equation}\label{eq:barotropic2} 
\left\{ 
\begin{aligned} 
&\partial_t\rho +\text{div}(\rho v) = 0, 
\\ 
&\partial_t v + \nabla_v v + \nabla\Big(2V - 2f(\rho) - \frac{2 \Delta\sqrt{\rho}}{\sqrt{\rho}} \Big) = 0 
\end{aligned} \right. 
\end{equation} 
for some functions $V\colon M\to \mathbb R$ and $f\colon (0,\infty) \to \mathbb R$. 
Then, the (time-dependent) complex-valued wave function $\psi=\sqrt{\rho e^{\mathrm{i}\theta}}$ 
given by the Madelung transform satisfies the nonlinear Schr\"{o}dinger equation on $M$: 
\begin{equation}\label{eq:schrodinger} 
\mathrm{i}\partial_t\psi = - \Delta\psi +  V\psi - f(|\psi|^2)\psi. 
\end{equation} 
The 1D Madelung transform, when interpreted in terms of curvature and torsion of the curve $\gamma$, 
reduces to the Hasimoto transform. 

It turns out that the Madelung transform not only maps one Hamiltonian equation to another, 
but it also preserves the symplectic structures related to the equations \cite{Ren12, KMM18}. 
More precisely, 
let $PC^\infty(M,\C)$ denote the complex projective space of smooth complex-valued functions $\psi$ on $M$:
its elements are cosets $[\psi]$ of the unit $L^2$-sphere of wave functions. 
The Madelung transform induces a symplectomorphism between 
$PC^\infty(M,\C\backslash \{0\})$, 
the projective space of non-vanishing complex functions,
and 
the cotangent bundle of probability densities $T^*{\rm Dens}(M)$ 
equipped with the canonical symplectic structure \cite{KMM18}. 
Furthermore, 
the Madelung transform is an isometry and a K\"ahler map between the spaces 
$T^\ast {\rm Dens}(M)$ 
equipped with the Sasaki-Fisher-Rao metric, which is the cotangent lift of the Fisher-Rao metric 
on the space of densities ${\rm Dens}(M)$, 
and 
$P C^{\infty}(M,\C\backslash \{0\})$ 
equipped with the Fubini-Study metric and the natural symplectic structures defined above, 
see \cite{KMM18, KMM20}. 
Finally, in \cite{Fu17} it was shown that the Madelung transform can be regarded as the momentum mapping 
for the space of wave functions regarded as half-densities on $M$ and acted upon by the semi-direct product group 
of diffeomorphisms and smooth functions. 
\begin{problem} 
Extend the above results on symplectomorphism and the K\"ahler map to the wave functions with zeros on $M$. 
Explain the quantization condition controversy {\rm \cite{Wallstrom, Fritsche}} in the language of the momentum map 
for the above semi-direct product group. 
\end{problem} 

The connection between equations of quantum mechanics
and hydrodynamics might shed some light on the hydrodynamic quantum analogues studied,
 e.g., in \cite{Couder, Bush}: 
 the motion of bouncing droplets in certain vibrating liquids manifests many properties of quantum mechanical particles. 
\begin{problem} 
Explain the droplet-quantum particle correspondence by a combination of  averaging and the Madelung transform. 
\end{problem}

\bigskip

\section{Mechanical models of direct and inverse cascades}
\subsection{Three models of energy propagation} 
\begin{figure}
\centering
\includegraphics[width=2.2in]{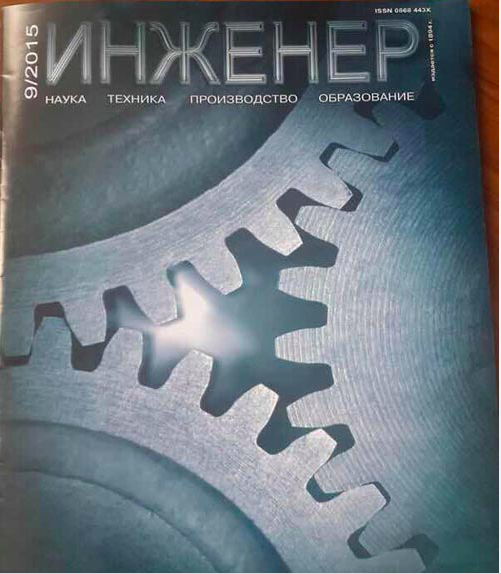}
\caption{\small  The 2D direct and inverse cascades together. (Illustration from the cover of Russian magazine ``Engineer", where the gears are Education, Science, and Enterprise. No comment.)}
\label{Fig.3}
\end{figure}

The inverse cascade phenomenon looks especially striking if we think of its mechanical model. 
Imagine a mechanism consisting of a countable number of wheels connected with gears, chains, springs 
and other joints which are assumed to be weightless and frictionless. 
Suppose that at the moment $t=0$ some wheels are set into rotation (Figure~\ref{Fig.3}).
In the course of motion the energy is redistributed among the wheels. The standard idea of statistical mechanics 
is that the energy tends to the uniform distribution between the wheels so that it will spread further and further. 
However, there are some other types of behaviour of such infinite mechanisms.

\medskip

On the one hand, the energy can spread so fast that at least part of it escapes to infinity in finite time 
and the total energy in the system decreases. 
On the other hand, the energy may become ``trapped'', i.e. it does not spread at all, 
and, moreover, it is concentrated in the first few wheels and its distribution does not depend on the initial energy profile, 
provided the energy is initially contained in any {\it finite number} of wheels.

There is also a softer regime, where the energy does not spread to all the wheels but its distribution depends 
on the initial profile (the Fermi--Pasta--Ulam regime). This might seem implausible but the fluid presents us 
with examples of such behaviour. 
In fact, consider a fluid flow $u(x,t)$ on the torus $\mathbb T^n$, $n=2$ or $n=3$. 
We can regard the Fourier coefficients $u_k(t)$ as the analogue of the $k$th wheel angular velocity. 
The time evolution of $u_k$ is described by a certain bilinear system of ODEs which can be regarded as 
a description of connections between the wheels (we could even design a realistically looking ``mechanism'' 
made of weightless and frictionless parts realizing these connections). 
Then, for $n=3$, we can expect a breakdown of a regular solution of the Euler equation 
and a transition to a turbulent motion with decreasing energy (which will be discussed below in more detail). 
If $n=2$, we observe the inverse energy cascade with the formation of a few large vortices 
so that  the energy is concentrated in a few lower harmonics, while the energy spectrum is decreasing over frequencies. 

\medskip

Thus, we have at least three types of behaviour of an ``infinite mechanism'':  
the tendency to the energy equidistribution, 
the inverse cascade (the energy tends to concentrate in the first few modes) 
and the direct cascade (the energy escapes to infinity in finite time). 
It is important to find out which properties of the mechanisms are responsible for so different behaviour. 
As a first step in this direction we can try to design some models, i.e. some simpler devices displaying similar behaviour. 
The original mechanism, such as a fluid in the Fourier representation, is too complicated to yield to the statistical theory 
with a lot of unrelated features. It would be interesting to find simpler (though infinite) mechanical systems 
which display the same statistical behaviour and which can be regarded as models of the fluid in this respect. 

\medskip

\subsection{A model of energy equidistribution} 
A classical example of a mechanical system with a well established tendency to the energy equidistribution 
is the gas of a large (but finite) number of solid balls moving inside a bounded domain.
It appears possible to modify this system to obtain a system with the direct energy cascade. 
Consider a system of a countable number of balls $B_j$ moving inside a bounded domain (``box'') $D$. 
Suppose these balls fall into a countable number of ``families''  $F_1, F_2, \ldots$ 
such that the family $F_i$ includes $n_i$ equal balls of radius $r_i$ and mass $m_i$. 
For simplicity suppose that the balls of a family $F_i$ ``feel'' only balls of the neighbouring families 
$F_{i-1}$ and $F_{i+1}$ and can penetrate through the balls from other families without any resistance. 
This means that the Hamiltonian of the system has the form 
$$ 
H(p,q) 
= 
\sum\limits_{i}\sum\limits_{F_i}\frac{p_j^2}{2 m_i} + \sum\limits_{i}\sum\limits_{\begin{subarray}{l}{B_j\in F_i} 
\\ 
{B_k\in F_{i+1}}\end{subarray}} U_{i, i+1}(|q_j - q_k|). 
$$
Suppose that the total mass is $M=\sum_i n_i m_i <\infty$ 
and the sequences of masses $m_i$ and radii $r_i$ are decreasing fast enough. 
Then each ball $B_j\in F_i$ is moving through a ``gas'' formed by the balls $B_k\in F_{i+1}$. 
The ``gas'' particles $B_k \in F_{i+1}$ are feeling the resistance of the ``gas'' formed by the balls $B_k\in F_{i+2}$, etc. 
The question is whether one can define the sequences $n_i, m_i, r_i$ in such a way that the direct energy cascade 
in the direction of growing $i$ would occur, and the energy would dissipate from the system? 
If the answer to the first question is affirmative, could one point at the details of an actual 3D fluid 
which play the role of the balls of different families? 

Now, the next natural question concerns what happens with the energy that has escaped from the system? 
One can introduce the ``limit absorption principle'' for our system. Namely, one can introduce a friction for the balls 
of the family $F_N$ which absorbs the energy and then let $N$ go to infinity. But this solution is not completely legitimate 
as here the energy sink is included from the beginning and it is only moved farther and farther away. 
So, this explanation of the energy escape is circular. 
\begin{problem} 
Is there a more natural way to introduce energy dissipation without explicit introduction of a friction mechanism?
\end{problem} 

\medskip

It would be intriguing to find a relation of the above model to  the mechanism of the energy dissipation in the weak solutions of the Euler equations whose rate is defined by the Duchon--Robert formula \cite{Duchon}. This would follow the footsteps of Maxwell's molecular vortex model for electromagnetic waves \cite{Max}.

\bigskip 

\begin{figure}
\centering
\includegraphics[width=5.6in]{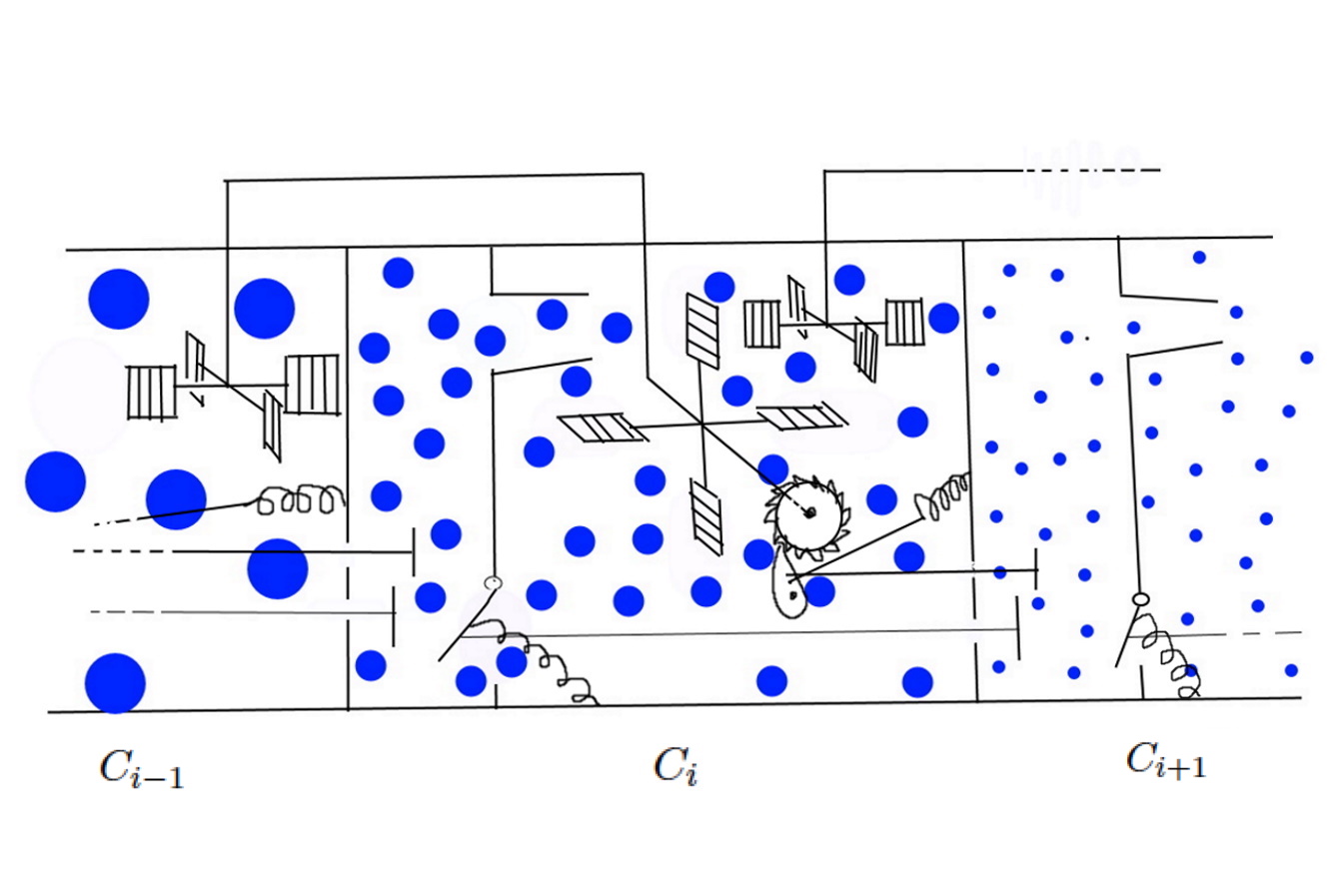}
\vspace{-17pt}
\caption{\small  A mechanical apparatus to mimic the inverse cascade.}
\vspace{-10pt}
\label{Fig.4}
\end{figure}
\medskip

\subsection{A model of the inverse cascade} 
A mechanical model of the inverse cascade requires a more sophisticated design.
First of all, the system will consist of  a countable number of chambers $C_1, C_2, \ldots$ with rigid walls, 
see Figure~\ref{Fig.4}.
Each chamber $C_i$ is separated by a vertical wall into two parts, call them the left room and the right room. 
The wall contains two doors, the upper and the lower door. Each chamber $C_i$ contains $n_i$ equal elastic balls $B_j^i$ 
of radius $r_i$ and mass $m_i$ moving inside $C_i$ and interacting with the walls, with the details of the mechanism 
and with other balls 
according to the laws of elastic collision. There is a shutter  at the lower door equipped with a spring 
which, when open, permits the balls to enter from the right room  into the left one
and, when shut, prevents the balls from going back. 
The shutter  is connected to a damper which is interacting with the balls in the next chamber $C_{i+1}$. 
The upper door  has no shutter and the balls can move freely through this door in both directions.

The balls enter the left room from the right one through the lower and the upper doors  and exit 
only through the upper door, provided the shutter at the lower door  works properly (this is, of course, a true Maxwell's demon). 
The last condition can be satisfied provided that the shutter, being a part of the system, is permanently cooled, 
i.e. its energy being transferred to the balls in the chamber $C_{i+1}$ 
(see the analysis of  Feynman \cite{Feynman}). 
To this end, the balls in $C_{i+1}$ should be much smaller and much more numerous than in $C_i$, 
i.e. $m_{i+1}\ll m_i$, $r_{i+1}\ll r_i$ and $n_{i+1}\gg n_i$. Then the balls in $C_{i+1}$ form a ``gas'' 
which is effectively viscous and absorbs the energy of the shutter. 
If  ``Maxwell's demon'' works properly then the balls in $C_i$ are entrained into the circular motion: 
on average, they enter from the right room into the left room through the lower door 
and leave from the right room mostly through the upper door. 
Thus, there appears a stream of balls from the upper door. 
Let us put a ``turbine''  which is rotated by this (possibly weak) stream. In order to ensure its rotation we put a ratchet-and-pawl which would prevent the reverse rotation of the turbine. To make this device work we should attach it to the second damper using the balls in $C_{i+1}$ to dissipate the energy.
Let this turbine drive through a system of connecting parts a ``stirrer''  
which transfers energy from the turbine  to the balls in the chamber $C_{i-1}$; 
these balls should be much larger than the ones in $C_i$, i.e. $m_{i-1}\gg m_i$, $r_{i-1}\gg r_i$, and $n_{i-1}\ll n_i$. 

The system should work as follows. 
The balls in the chamber $C_i$ are, on the average, taking part in the circular motion entering 
from the right room  into the left one through the lower door and from the left room into the right one through the upper door (equipped with a nozzle), 
while the shutter at the lower door is damped by the damper. The latter is braked by the ``gas'' of the balls in $C_{i+1}$ 
which are much smaller than the balls in $C_i$. 
The energy of the stream of the balls in $C_i$ is transferred through the turbine 
to the much larger balls in the chamber $C_{i-1}$. 
The ratchet-and-pawl pair is cooled by a similar damper (see the above analysis of this pair by Feynman).
As a result, on the average the energy is transferred to the balls in the first few chambers. 

This mechanism looks like a sort of {\it perpetuum mobile}. 
However, it is neither a perpetuum mobile of the first nor of the second kind. 
In fact, it is not a perpetuum mobile at all 
but, rather, it is a chain of heat engines: 
the balls inside each chamber $C_i$ are the ``working body'' of the engine. 
The balls in the next chamber $C_{i+1}$ play the role of  a cooler (they are cooling the shutter, the  
``Maxwell's demon'' and the ratchet-and-pawl), 
while  the balls in the previous chamber $C_{i-1}$ are playing the role of the load. 
The engines are quite primitive and their efficiency is very low; 
however, there is an infinite number of them so that their overall efficiency is 100\%. 
Hence, here is the problem. 
\begin{problem} 
Is there a similar mechanism which describes a cascade in real hydrodynamics? 
Is it possible to define the parameters $n_i, \, m_i,$ and $ r_i$ in such a way that the above apparatus works as intended? 
\end{problem} 

It would be interesting and important to find a link between this device and a more convenient heat engine; 
in particular, to find some analogues of Maxwell's demon 
or, perhaps even more importantly, 
to show that a similar mechanism works in a 2D ideal fluid, thus ensuring inverse energy cascade in it. 
\bigskip

{\bf Acknowledgments.} We are indebted to the organizers of the IAS Symplectic Dynamics Program, where this paper was conceived, and in particular to Helmut Hofer for his constant encouragement.  We are also grateful to Theodore Drivas and to the anonymous referee for many useful suggestions.
B.K. was partially supported by a Simons Fellowship and an NSERC Discovery Grant. G.M. gratefully acknowledges support from the SCGP, Stony Brook University, where part of the work on this paper was done. 
Data sharing is not applicable to this article as no datasets were generated or analyzed during the current study.


\bigskip



\end{document}